\documentclass[twocolumn,10pt]{asme2ej}

\usepackage{amsmath}
\usepackage{amsfonts}
\usepackage{bm}
\usepackage{dsfont}
\usepackage{enumitem}
\usepackage{graphicx}
\usepackage{float}
\usepackage{rotating}
\usepackage{booktabs}
\usepackage{multirow}
\usepackage{xcolor}

\title{Robust reliability-based topology optimization under random-field material model}

\author{Trung Pham\thanks{Address all correspondence related to this paper to this author.}
    \affiliation{
    	Department of Aerospace Engineering \\
    	University of Michigan \\
    	Ann Arbor, Michigan 48109, USA \\
        Email: trungp@umich.edu
    }	
}

\author{Christopher Hoyle
    \affiliation{ 
        School of Mechanical, Industrial \& Manufacturing Engineering \\
    	Oregon State University \\
    	Corvallis, Oregon 97331--6001, USA \\
        Email: chris.hoyle@oregonstate.edu
    }
}

\begin{document}

\maketitle

\begin{abstract}
{\it This paper proposes an algorithm to find robust reliability-based topology optimized designs under a random-field material model. The initial design domain is made of linear elastic material whose property, i.e., Young's modulus, is modeled by a random field. To facilitate computation, the Karhunen$-–$Lo\`{e}ve expansion discretizes the modeling random field into a small number of random variables. Robustness is achieved by optimizing a weighted sum of mean and standard deviation of a quantity of interest, while reliability is employed through a probabilistic constraint. The Smolyak-type sparse grid and the stochastic response surface method are applied to reduce computational cost. Furthermore, an efficient inverse-reliability algorithm is utilized to decouple the double-loop structure of reliability analysis. The proposed algorithm is tested on two common benchmark problems in literature. Finally, Monte Carlo simulation is used to validate the claimed robustness and reliability of optimized designs.}
\end{abstract}

\section{Introduction}
\label{intro}
A mechanical structure is characterized by its boundary and loading conditions, its material properties, and its topology. Finding an appropriate topology is often a major task in structural design, which has fueled the rise of topology optimization (TO) in the last two decades. Without considering a designer's experience, TO is a mathematical tool to identify the optimal size, shape, and connectivity of a design \cite{bendsoe_topology_2003}, resulting in improved performance while using the least amount of material. However, research in TO often only concerns with deterministic inputs, while uncertainty is inherent in nature, which manifests itself in the stochasticity of random parameters of engineered systems. The analysis and design of engineered systems are affected heavily by uncertainty; for example, modern design codes, such as ACI 318 \cite{aci318} and AISC 360 \cite{aisc360}, have comprehensive recommendations of safety factors for loading, material property, construction conditions, etc., which obviously are intended to take into account uncertainty. Among different sources of uncertainty, material property is intrinsically random in space, which has been modeled by random field in the design of composite structures \cite{sriramula_experimental_2013,nader_probabilistic_2008,sriramula_quantification_2009}. Such a modeling technique has been especially popular in the vast literature of the Stochastic Finite Element Method \cite{sudret_2000, ghanem2003stochastic,charmpis_need_2007}, which certainly proves its validity and the need to be considered in TO. So far, uncertainty in TO has been treated separately by robust optimization or reliability-based optimization, while both robustness and reliability are desired properties of design under uncertainty. Therefore, this paper presents an algorithm to find robust reliability-based topology optimized design under a random-field material model.

There are a number of steps in our proposed algorithm, which are addressed in depth in the subsequent sections. Here we provide a brief overview of them. First, from a known covariance function, the modeling random field is estimated by a random polynomial using the Karhunen$-–$Lo\`{e}ve expansion. To make the design robust, a weighted sum of mean and standard deviation of a quantity of interest, which are computed by a Smolyak-type sparse grid, is considered as the objective function. Reliability of the design is reflected in the probabilistic constraint, which is handled by the Sequential Optimization and Reliability Assessment (SORA) method \cite{du_sequential_2004}$-$a single-loop inverse-reliability algorithm$-$coupled with the performance measure approach \cite{tu_new_1999} to reduce the computational cost of reliability analysis. The stochastic response surface method approximates the random output, which is required to solve the inverse reliability analysis problem.

The layout of this paper is as follows. Section \ref{background} is a literature review of uncertainty propagation, robust and reliability-based optimization, and reliability analysis, with focus on TO under uncertainty. Section \ref{tou} derives the mathematical formulation of the deterministic TO problem and the robust reliability-based topology optimization (RRBTO) problem. This section also exposes the details of our proposed solution for the RRBTO problem. Two numerical examples show how our proposed algorithm works, and are verified by Monte Carlo simulation in section \ref{results} followed by discussions of the results. Finally, the paper is summarized with major findings, and then suggests future work.

\section{Background}
\label{background}
The earliest idea of topology optimization (TO) can be traced back to Michell's paper \cite{michell1904} in 1904. Since then TO has been mature enough to have its own treatise \cite{bendsoe_topology_2003}. As a mathematical optimization problem, TO requires specification of objective function(s) and constraint(s), which do not involve any probabilistic quantities when using deterministic inputs. Thus, changes are needed to deal with uncertainty in the forms of robust optimization (RO) and reliability-based optimization (RBO). A number of papers, which are reviewed below, have tried to integrate uncertainty into a TO problem using RO and RBO separately.

RO primarily aims to minimize the variability of an output of interest \cite{taguchi_1986}, due to uncertainty, around its mean value. Therefore, this goal can be formulated by minimizing a weighted sum of mean and standard deviation of the output of interest. This approach was chosen in several papers covering various sources of uncertainty and solution methods: spatial variation of manufacturing error with Monte Carlo simulations \cite{schevenels_robust_2011}; random-field truss material with a multi-objective approach \cite{richardson_robust_2015}; random loading field and random material field with the level set method \cite{chen_level_2010}; material and geometric uncertainties with stochastic collocation methods and perturbation techniques \cite{lazarov_topology_2012,lazarov_topology_2012-1}; misplacement of material and imperfect geometry \cite{jansen_robust_2013,jansen_robust_2015}; Young's modulus of truss members with a perturbation method \cite{asadpoure_robust_2011}; random-field material properties with a polynomial chaos expansion \cite{tootkaboni_topology_2012}; geometric and  material property uncertainties with a stochastic perturbation method for frame structures \cite{changizi_robust_2017}; material uncertainty with known second-order statistics \cite{da_silva_topology_2016}; random spatial distribution of Young's modulus and loading uncertainty in a stress-based problem \cite{da_silva_stress-based_2017,dasilva_topology_2018}; and random loading field with stochastic collocation methods \cite{zhao_robust_2015}. The robust topology optimization (RTO) problem is solved by a unified framework based on polynomial chaos expansion in \cite{richardson_unified_2016}, while \cite{zhao_robust_2014} tackled the problem exploiting the linear elasticity of structure. The seemingly arbitrary factors in the weighted sum are often cited as one major weakness of this RO methodology \cite{jalalpour_efficient_2016,richardson_robust_2015}; however, they are well-defined in decision-based design reflecting risk-taking attitude of designers \cite{hoyle_2014,lewis_2006,beck_comparison_2015}.

Instead of modifying the objective function, RBO makes some of the constraints probabilistic$-$probability of failure event is used in place of the event itself. This change requires specialized methods to handle, because the probabilistic constraints are expressed by multiple integrals of the joint probability density function (PDF) of random variables, both of which are either practically impossible to obtain or very difficult to evaluate \cite{achintya_2000}. Many methods have been devised to overcome such difficulties, which were surveyed thoroughly in \cite{valdebenito_survey_2010}. Within the scope of this paper, we only briefly review the first-order reliability methods (FORM), the second-order reliability methods (SORM), and the Sequential Optimization and Reliability Assessment (SORA) method. FORM appeared early \cite{cornell_1969} together with the concept of reliability index \cite{hasofer_1974} to solve RBO problems. SORM \cite{fiessler_1979} followed to improve accuracy of the FORM in case of highly nonlinear limit state functions and/or slow decay of the joint PDF. The main idea of FORM and SORM is to approximate the limit state functions using first-order and second-order Taylor series, respectively, at appropriate values (i.e., means) of random variables. This results in a double-loop optimization problem to find the most probable point (MPP). In the context of reliability-based topology optimization (RBTO), directly solving the double-loop optimization problem has been shown in \cite{maute_reliability-based_2003} for MEMS mechanisms with stochastic loading, boundary conditions as well as material properties; in \cite{sato_reliability-based_2018} for shape uncertainty; in \cite{kang_reliability-based_2018} for geometric imperfections; in \cite{mogami_reliability-based_2006} for frame structures using system reliability under random-variable inputs; in \cite{kang_reliabilitybased_2004} for electromagnetic systems; in \cite{jung_reliability-based_2004} for geometrically nonlinear structures; in \cite{luo_reliability_2014} for local failure constraints; and in \cite{daSilva2018} for continuum structures subject to local stress constraints. In \cite{papadimitriou_reliability-based_2018}, FORM was replaced by a mean-value, second-order saddlepoint approximation method, which was asserted to be more accurate. The double-loop approach is prohibitively expensive and lacks robustness when a large number of random variables presents \cite{schueller_critical_2004}. For this reason, single-loop approaches have been developed, in which, i.e., the Karush$-$Kuhn$-$Tucker (KKT) optimality conditions are utilized to avoid the inner loop. Both \cite{nguyen_single-loop_2011} and \cite{silva_component_2010} used variants of the single-loop method in \cite{liang_single-loop_2004} for component and system reliability-based TO. In \cite{kharmanda_reliability-based_2004}, it is somewhat unique when the authors used their own single-loop method \cite{kharmanda_efficient_2002}. Kogiso et al. \cite{kogiso_reliability-based_2010} applied the single-loop-single-vector method \cite{chen_reliability_1997} for frame structures under random-variable loads and nonstructural mass. Another way to bypass the double-loop problem is the decoupling approaches \cite{valdebenito_survey_2010}, in which reliability analysis results are used to facilitate the optimization loops. Among them, the SORA method is known for its simple implementation compared to the above single-loop methods, and its efficiency with FORM \cite{du_sequential_2004,lopez_reliability-based_2012}. This method was employed for RBTO under random-variable inputs in \cite{zhao_reliability-based_2015} and \cite{zhao_comparison_2016}. Meta modeling or surrogate modeling used together with simulation techniques to solve RBO problems has received considerable attention \cite{wang_review_2007}, but still remained relatively unexplored in TO literature. In \cite{patel_classification_2012}, reliability was assessed using a probabilistic neural network classifier for truss structures under random Young's modulus. 

Both robustness and reliability are desired properties of design under uncertainty; however, to the best of our knowledge, this paper is the first one considering both criteria in TO. Therefore, a literature review of robust reliability-based optimization (RRBO) has to be drawn from other fields. The RRBO problem was investigated in \cite{du_integrated_2004} using an inverse reliability strategy; in \cite{youn_performance_2005} using a performance moment integration method to estimate the product quality loss; in \cite{mourelatos_methodology_2006} using a preference aggregation method to produce a single-objective RBO problem; and in \cite{tang_sequential_2012} under epistemic uncertainty. Both \cite{forouzandeh_shahraki_reliability-based_2014} and \cite{rathod_optimizing_2013} used a genetic algorithm to solve the problem. The dimension reduction method and its derivatives were introduced in \cite{lee_dimension_2008,motasoares_stochastic_2006,youn_reliability-based_2009} as an alternative approach to the RRBO problem.

With respect to RBTO, the approaches in \cite{zhao_reliability-based_2015}, \cite{jalalpour_efficient_2016}, and \cite{keshavarzzadeh_topology_2017} are closest to ours. Still, random-field modeling was not used for material property in \cite{zhao_reliability-based_2015} and \cite{keshavarzzadeh_topology_2017}. Furthermore, several concerns can be identified from \cite{keshavarzzadeh_topology_2017}. One of the most important stages in their method is the approximation of both failure probability and its sensitivity, which were calculated by Monte Carlo sampling. Direct Monte Carlo sampling is well-known to have variability \cite{taflanidis_efficient_2008}, meaning two independent runs are very likely to get different values of failure probability and its sensitivity which would obviously affect the optimization results. Another concern is that the value of the parameter $\epsilon$ needed to replace the Heaviside function with a smooth approximation \cite{keshavarzzadeh_gradient_2016} and was chosen by a ``recommendation" backed by observation only. In \cite{jalalpour_efficient_2016}, the modeling random field was assumed with a known marginal distribution, and in order to apply a perturbation technique, random variability of Young's modulus had to be small. Both of these assumptions clearly restrict the general applicability of their method. Lastly, robustness against uncertainty was not studied in the three papers. As described in the following sections, our proposed method considers random field uncertainty with the Karhunen$-–$Lo\`{e}ve (KL) expansion used to reduce the dimension of the random field. The KL expansion covers a large class of random field without any restrictions on random variability. In this way, we are able to use the FORM-based inverse reliability method within the SORA framework coupled with the stochastic response surface method to avoid the aforementioned weaknesses of direct Monte Carlo sampling.

\section{Topology Optimization under Uncertainty}
\label{tou}
\subsection{Deterministic Topology Optimization}
\label{deter}
A standard notation is adopted throughout this paper$-$bold upper and lower case letters denote matrices and vectors, respectively. The below formulation shows a density-based deterministic topology optimization (DTO):
\begin{equation}
    \begin{aligned}
        \min_{\boldsymbol{\rho}} && &C( \boldsymbol{\rho} ) = \text{\bf{u}}^T \text{\bf{K}} \text{\bf{u}}\\
        \text{subject to} && &\text{\bf{K}}(\boldsymbol{\rho}) \text{\bf{u}}(\boldsymbol{\rho}) = \text{\bf{f}}, \\
        && & \frac{V(\boldsymbol{\rho})}{V_0} = \gamma, \\
        && &0 < \rho_{\text{min}} \leq \boldsymbol{\rho} \leq 1,
        \label{e:1}
    \end{aligned}
\end{equation}
where $\boldsymbol{\rho}$ is the vector of design variables of the TO problem, which also are the deterministic finite-element densities; $V( \boldsymbol{\rho} )$ and $V_0$ are the total and the initial volume of the finite-element mesh, respectively; $\gamma$ is the predetermined volume fraction; $\text{\bf{K}}(\boldsymbol{\rho})$, $\text{\bf{u}}(\boldsymbol{\rho})$, and $\text{\bf{f}}$ are the stiffness matrix, the displacement vector, and the external load vector, respectively; and $C( \boldsymbol{\rho} )$ is the structural compliance. In the optimization problem~(\ref{e:1}), there are three constraints: the first expresses the equilibrium of the structure; the second requires the optimized design to have a prescribed volume; and the third is a component-wise inequality, in which each density (design variable) must be between 1 and a lower limit (i.e., $\rho_{\text{min}} = 0.001$).

To ensure manufacturability, the optimized design must have a well-defined boundary, which is not guaranteed if solving (\ref{e:1}) directly because there is nothing to prevent intermediate values of densities from dominating the design. Hence, the Solid Isotropic Material with Penalization (SIMP) method \cite{bendsoe_topology_2003} is used to make intermediate densities unfavorable compared to $\rho_{\text{min}}$ or 1. According to SIMP, each finite element has a Young's modulus $E_i$ specified by $E_i = \rho_i^pE_i^0$, where $p$ is the penalization factor and $E_i^0$ is the initial value of the Young's modulus corresponding to unit density. The interpretation and possible values of $p$ were elaborated in \cite{bendsoe_material_1999}. Any established gradient-based algorithms can solve the problem~(\ref{e:1}) after it is converted into a nonlinear optimization problem using the SIMP method. This paper follows common practice in the literature, selecting the Method of Moving Asymptotes (MMA) \cite{svanberg_method_1987,svanberg_class_2002} as the optimizer of the DTO problem. The MMA has proved its reliability and competitive performance in various settings of TO. However, SIMP alone is plagued with checkerboarding, mesh dependence, and local minima \cite{sigmund_numerical_1998}. Many mesh-independent filtering methods \cite{sigmund_morphology-based_2007} have been designed to preclude checkerboarding and mesh dependence, while local minima remain an open question. This paper uses the density filtering \cite{bruns_topology_2001,bourdin_filters_2001} as implemented in \cite{andreassen_efficient_2011}. The next sections describe how uncertainty shapes our problem formulation and the solution algorithm.

\subsection{Robust Reliability-based Topology Optimization}
\subsubsection{Problem Formulation}
Considering input uncertainty modeled by a random field $y(\omega,\mathbf{x})$, a robust reliability-based topology optimization (RRBTO) problem is formulated as follows:
\begin{equation}
    \begin{aligned}
        \min_{\boldsymbol{\rho}} && &\kappa_1 \mu\left[C( \boldsymbol{\rho},y )\right] + \kappa_2 \sigma\left[C( \boldsymbol{\rho},y )\right] \\
        \text{s.t.} && &\text{\bf{K}}(\boldsymbol{\rho},y) \text{\bf{u}}(\boldsymbol{\rho},y) = \text{\bf{f}}, \\
        && & \frac{V(\boldsymbol{\rho})}{V_0} = \gamma, \\
        && &P_i\left[g_i( \boldsymbol{\rho},y) < 0\right] \leq P_i^0, \; i = 1, 2\dots, m,\\
        && &0 < \rho_{\text{min}} \leq \boldsymbol{\rho} \leq 1,
        \label{e:2}
    \end{aligned}
\end{equation}
where $\mathbf{x} \in D \subset \mathds{R}^d$ is coordinates of a point in a $d$-dimensional physical domain $D$; $\omega \in \Omega$ is an element of the sample space $\Omega$; $\mu\left[C( \boldsymbol{\rho},y )\right]$ and $\sigma\left[C( \boldsymbol{\rho},y )\right]$ are the mean and standard deviation of the compliance $C( \boldsymbol{\rho},y )$, respectively; $\kappa_1$ and $\kappa_2$ are the real non-negative weighting factors. The limit state function $g_i(\boldsymbol{\rho},y)$ is defined so that $g_i(\boldsymbol{\rho},y) < 0$ means failure of the design, and $P_i[g_i(\boldsymbol{\rho},y) < 0]$ shows the probability of the $i^{th}$ failure event. The target probability $P_i^0$ is the upper bound of the failure probability $P_i$ and often defined as $P_i^0=\Phi(-\beta_i)$, where $\beta_i$ is the reliability index and $\Phi(\cdot)$ is the standard normal cumulative distribution function. In this paper the random field $y(\omega,\mathbf{x})$ is taken to be the material Young's modulus, which must be physically meaningful (i.e., taking only positive values) and is modeled as in \cite{lazarov_topology_2012}:
\begin{equation}
    E(\mathbf{x}) = F^{-1} \circ \Phi\left[y(\omega,\mathbf{x})\right],
    \label{e:young}
\end{equation}
where $\Phi[\cdot]$ is the standard normal cumulative distribution function (CDF), and $F^{-1}$ is the inverse of a prescribed CDF. The uniform distribution is chosen for the two numerical examples resulting in
\begin{equation}
    E(\mathbf{x}) = a + (b-a)\Phi\left[y(\omega,\mathbf{x})\right],
    \label{e:uniform}
\end{equation}
where $a$ and $b$ are the two bounds of the distribution. The log-normal and the beta distribution are also capable of modeling non-negative, bounded physical quantities, which can supersede the uniform distribution in~(\ref{e:young}) with minimal effort.

A number of challenges need to be cleared before we are able to solve~(\ref{e:2}). The modeling random field, which tries to capture spatial variability of material property, needs to be cast into an explicit, computable form because a defined value of material property is required to perform finite element analysis. The Karhunen$-–$Lo\`{e}ve (KL) expansion in Section~\ref{kle} is able to turn a random field into a series of random variables. The mean and standard deviation of the compliance, and the probabilistic constraints are often very hard or expensive to evaluate because of complex geometry of their domains. A Smolyak-type sparse grid in Section~\ref{ssg} is an efficient method to calculate the mean and standard deviation, while Inverse Reliability Analysis (IRA) and SORA in Section~\ref{irs}, coupled with the Stochastic Response Surface Method (SRSM) in Section~\ref{srsm}, handle the probabilistic constraints effectively by avoiding the demanding double-loop problem. Combining the above methods, a detailed description of our proposed algorithm is laid out in Section~\ref{algorithm}.

\subsubsection{Smolyak-type Sparse Grid}
\label{ssg}
In practice, it is very difficult or even impossible to calculate the mean and standard deviation of the compliance in (\ref{e:2}) analytically through multidimensional integrals. Such difficulties have motivated the development of various numerical methods such as simulation-based methods (i.e., Monte Carlo (MC), important sampling, adaptive sampling, etc.), and the stochastic collocation methods (SCM) \cite{lee_comparative_2009}. Simulation-based methods are usually more straightforward to implement and embarrassingly parallel, and their cost does not depend on the number of dimensions; however, even with better sampling techniques, they still require a lot more sampling points than the SCM. Depending on the smoothness of the target function, the convergence rate of the SCM can be orders of magnitude faster than MC-based methods \cite{xiu2010numerical,xiong_new_2010}. The SCM approximates the quantity of interest by a weighted sum, whose weighing factors and terms are computed at specific collocation points. Locating such points is one of the central topics in SCM. The popular approach is to pick a known one-dimensional quadrature rule and then build up the multidimensional grid from the one-dimension rule. Interested readers can find in \cite{gerstner_numerical_1998} a list of popular quadrature rules. The Gauss-Hermite quadrature, which is particularly suitable for approximating the mean of a normal distribution, is selected in this paper. The multidimensional grid can be constructed using a tensor product; nevertheless, due to the well-known curse of dimensionality, the cost of SCM on a full tensor-product grid is still excessively high for a large number of dimensions. The Smolyak-type sparse grid (SSG), which may be traced back to the Smolyak algorithm \cite{smolyak_1963}, can significantly reduce the cost by using only a subset of the full tensor grid. In \cite{xiu_high-order_2005} numerical experiments with random input, whose dimension was up to 50, showed that the SCM on sparse grids was more efficient than MC.

For the sake of completeness the SSG is reviewed here.  The construction presented below follows \cite{maitre2010spectral}. Consider a $d$-dimensional function $f(\mathbf{x})$, the difference $\Delta_l^{(1)}f$ is defined as
\begin{equation}
	\Delta_l^{(1)}f = \left(Q^{(1)}_l-Q^{(1)}_{l-1}\right)f,
	\label{e:sparse_grid_1}
\end{equation}
where
\begin{equation}
	\begin{aligned}
		Q^{(1)}_l f & =\sum_{i=1}^{n_l^{(1)}}f\left(\mathbf{x}_l^{(i)}\right)w_l^{(i)}, \\
		Q_0^{(1)}f  & =0,                                                               
	\end{aligned}
\end{equation}
and $n_l^{(1)}$ is the number of nodes for the $l$-level quadrature formula. The collocation points $\mathbf{x}_l^{(i)}$ and the corresponding weights $w_l^{(i)}$ are calculated using the Gauss-Hermite quadrature rule, which is a natural choice for the class of integrals involving an exponential function over an infinite interval such as mean and standard deviation \cite{davis2007methods}. The sparse integration formula at level $l$ is expressed as
\begin{equation}
	\begin{aligned}
		Q_l^{(d)} f = \sum_{|\ell| \leq l+d-1}\left(\Delta_{l_1}^{(1)} \otimes \Delta_{l_2}^{(1)} \otimes \cdots \otimes \Delta_{l_d}^{(1)} \right)f, 
		\label{e:sparse_grid_2}                         
	\end{aligned}
\end{equation}
where $|\ell|=l_1+l_2+\ldots+l_d$. The mean and standard deviation of $f(\mathbf{x})$ can be approximated using (\ref{e:sparse_grid_2}).

\subsubsection{Karhunen$-–$Lo\`{e}ve Expansion}
\label{kle}
In a physical system, the quantity of interest can be measured at spatial points over the system domain. If a random field is considered an appropriate model for such a quantity, it is then required to construct the random field from measurements. Several methods, including the Expansion Optimal Linear Estimator \cite{li_optimal_1993} and polynomial chaos expansion \cite{xiu2010numerical,ghanem2003stochastic}, have been adopted. Compared to others, the Karhunen$-–$Lo\`{e}ve (KL) expansion \cite{loeve2017probability} is ``the most efficient in terms of the number of random variables required for a given accuracy" \cite{sudret_2000}. The KL expansion of a random field $y(\omega,\mathbf{x})$ is given as
\begin{equation}
    y(\omega,\mathbf{x}) = E[\mathbf{x}] + \sum_{i=1}^\infty\sqrt{\lambda_i}\xi_i(\omega)e_i(\mathbf{x})
    \label{e:kl_2}
\end{equation}
where $E[\mathbf{x}]$ is the mean of the random field. The orthogonal eigenfunctions $e_i(\mathbf{x})$ and the corresponding eigenvalues $\lambda_i$ are solutions of the  following eigenvalue problem:
\begin{equation}
    \int_{D}K(\mathbf{x}_1,\mathbf{x}_2)e_i(\mathbf{x})d\mathbf{x}=\lambda_ie_i(\mathbf{x}) \quad \mathbf{x},\mathbf{x}_1,\mathbf{x}_2 \in D
    \label{e:kl_3}
\end{equation}
where $K(\mathbf{x}_1,\mathbf{x}_2)$ is the covariance function of the random field
\begin{equation}
    K(\mathbf{x}_1,\mathbf{x}_2) = E\left[y(\mathbf{x}_1)y(\mathbf{x}_2)\right] \quad \mathbf{x}_1,\mathbf{x}_2 \in D
\end{equation}
The random variables $\xi_i(\omega)$ are uncorrelated and satisfy:
\begin{equation}
    \begin{aligned}
        & E[\xi_i] = 0, E[\xi_i\xi_j] = \delta_{ij}, \\
        & \xi_i(\omega) = \frac{1}{\sqrt{\lambda_i}}\int_D\left(y(\omega,\mathbf{x})-E[\mathbf{x}]\right)e_i(\mathbf{x})d\mathbf{x},
    \end{aligned}
\end{equation}
where $\delta_{ij}$ is the  Kronecker delta. The infinite series in~(\ref{e:kl_2}) has to be truncated to use in practice. Because the influence of higher order terms decays rapidly, satisfactory precision can be achieved using only the first few terms of the expansion.

The KL expansion requires the solution of the eigenvalue problem~(\ref{e:kl_3}), which is pretty straightforward in the case of a random process (1-dimensional random field) \cite{ray_numerical_2013}. For the purpose of demonstration and without loss of generality, this paper in Section~\ref{results} assumes the separability of the covariance function of a 2-dimensional random field:
\begin{equation}
    \begin{aligned}
    K(\mathbf{s},\mathbf{t}) &= exp\left(\frac{-|s_1-t_1|}{l_1}\times\frac{-|s_2-t_2|}{l_2}\right) \\
    &=exp\left(-\frac{|s_1-t_1|}{l_1}\right)exp\left(-\frac{|s_2-t_2|}{l_2}\right), \\
    &\mathbf{s},\mathbf{t} \in D \subset \mathds{R}^2,
    \label{e:covariance}
    \end{aligned}
\end{equation}
where $l_1$ and $l_2$ are the correlation lengths in the two coordinate directions. The separability of the covariance function leads to separable eigenvalues and eigenfunctions, which are the product of their univariate counterparts \cite{wang_2008}.

\subsubsection{Inverse Reliability Analysis and SORA}
\label{irs}
\begin{figure}[ht]
    \centering
    \includegraphics[scale=0.85]{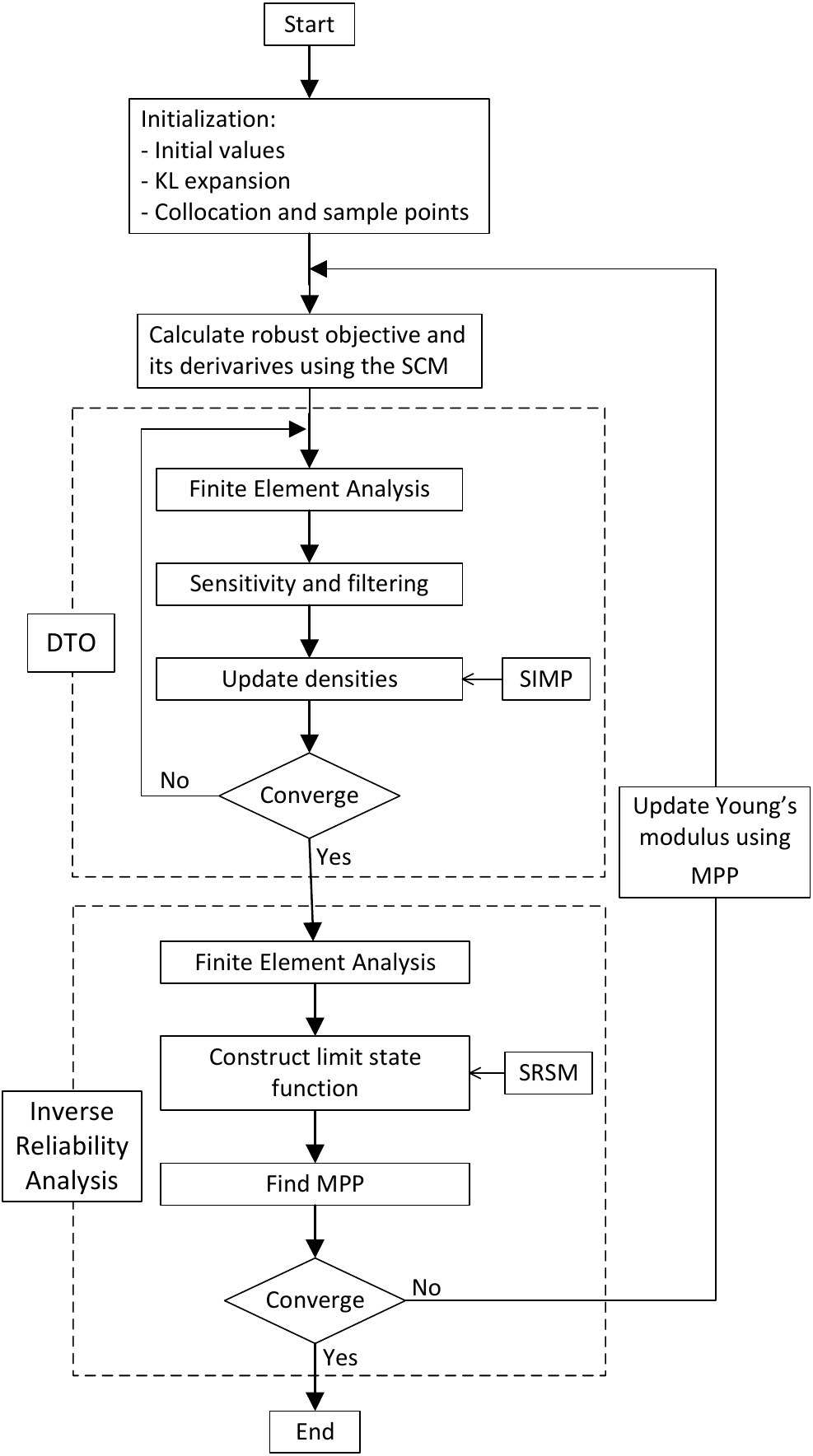}
    \caption{SORA-based RRBTO flowchart \cite{zhao_reliability-based_2015}.}
    \label{fig:rr_sora}
\end{figure}

The probabilistic constraints formulated as in~(\ref{e:2}) are called the reliability index approach (RIA); however, as \cite{tu_new_1999} reported, the performance measure approach (PMA) provides better numerical stability and higher rate of convergence. Using the PMA, (\ref{e:2}) is transformed as follows:
\begin{equation}
    \begin{aligned}
        \min_{\boldsymbol{\rho}} && &\kappa_1 \mu\left[C( \boldsymbol{\rho},y )\right] + \kappa_2 \sigma\left[C( \boldsymbol{\rho},y )\right] \\
        \text{s.t.} && &\text{\bf{K}}(\boldsymbol{\rho},y) \text{\bf{u}}(\boldsymbol{\rho},y) = \text{\bf{f}}, \\
        && & \frac{V(\boldsymbol{\rho})}{V_0} = \gamma, \\
        && &g_i( \boldsymbol{\rho},y) \geq 0, \; i = 1, 2\dots, m,\\
        && &0 < \rho_{\text{min}} \leq \boldsymbol{\rho} \leq 1.
        \label{e:ira_1}
    \end{aligned}
\end{equation}
The most notable change is that the probabilistic constraints are replaced by inequalities of the limit state functions. Solving the new problem requires a truncated KL expansion $y(\omega,\mathbf{x}) \approx y(\boldsymbol{\xi}(\omega),\mathbf{x})$, FORM, and inverse reliability analysis (IRA). In order to apply FORM, the random vector $\bm{\Xi}=\{\xi_i\}$ is first transformed into a vector of standard normal random variables $\bm{\Psi}=\{\psi_i\}$ using the Rosenblatt or the Nataf transformation $\bm\Psi = T(\bm{\Xi})$ or $\bm{\Xi} = T^{-1}(\bm{\Psi})$. Then, the most probable point (MPP) $\bm{\xi}_i^*$ in physical space or $\bm{\psi}_i^*$ in transformed space is obtained by solving the following IRA problem:
\begin{equation}
    \begin{aligned}
        \min_{\boldsymbol{\psi}} && &g_i( \boldsymbol{\psi} ) \\
        \text{s.t.} && &\parallel\boldsymbol{\psi}\parallel=\beta_i,
        \label{e:relia3}
    \end{aligned}
\end{equation}
where $g_i( \bm{\psi} )$ is the $i^{th}$ limit state function in transformed space. In this paper, the Matlab CODES toolbox \cite{codes_toolbox} is chosen to solve~(\ref{e:relia3}). Furthermore, the SORA framework is adopted to decouple the double-loop structure of~(\ref{e:ira_1}). In SORA, instead of nesting the optimization problem~(\ref{e:relia3}) within~(\ref{e:ira_1}), it serializes~(\ref{e:ira_1}) into a chain of loops of DTO and IRA (Fig.~\ref{fig:rr_sora}). Each $k^{th}$ loop starts with DTO followed by IRA:
\begin{equation}
    \begin{aligned}
        \min_{\boldsymbol{\rho}^k} && &\kappa_1 \mu\left[C( \boldsymbol{\rho}^k,y )\right] + \kappa_2 \sigma\left[C( \boldsymbol{\rho}^k,y )\right] \\ \\
        \text{s.t.} && &\text{\bf{K}}(\boldsymbol{\rho}^k,y) \text{\bf{u}}(\boldsymbol{\rho}^k,y) = \text{\bf{f}}, \\
        && & \frac{V(\boldsymbol{\rho}^k)}{V_0} = \gamma, \\
        && &g_i\left( \boldsymbol{\rho}^k,y(\boldsymbol{\xi}_i^{*(k-1)},\mathbf{x})\right) \geq 0, \; i = 1, 2\dots, m,\\
        && &0 < \rho_{\text{min}} \leq \boldsymbol{\rho}^k \leq 1.
        \label{e:relia4}
    \end{aligned}
\end{equation}
where $\boldsymbol{\xi}_i^{*(k-1)}$ denotes the MPP in physical space of $i^{th}$ limit state function in the $(k-1)^{th}$ loop. Solving~(\ref{e:relia4}) gives $\boldsymbol{\rho}^{*(k)}$, which is substituted into~(\ref{e:relia3}) to find the next MPP $\boldsymbol{\xi}_i^{*(k)}$ in the form of $\boldsymbol{\psi}_i^{*(k)}$:
\begin{equation}
    \begin{aligned}
        \min_{\boldsymbol{\psi}} && &g_i(\boldsymbol{\rho}^{*(k)},\boldsymbol{\psi} ) \\
        \text{s.t.} && &\parallel\boldsymbol{\psi}\parallel=\beta_i.
        \label{e:relia5}
    \end{aligned}
\end{equation}

The SORA framework coupled with the PMA can save computational cost significantly by reducing the number of reliability analyses performed to reach convergence of both DTO and IRA. However, as a heuristic method, the optimization solution $\boldsymbol{\rho}^*$ acquired by SORA may not be the one found in the double-loop problem, which is corrected in each loop by shifting the random design variables using information from the previous loop~\cite{yin_enhanced_2006}. Such modification is not necessary in this paper because there are only deterministic design variables (only the material property is random).

\subsubsection{Stochastic Response Surface Method}
\label{srsm}
As uncertainty is propagated from random input to output through complex computations such as finite element analysis, it is almost unlikely to derive the output directly as an explicit expression of input. However, such an expression is needed to run the gradient-based algorithm in IRA. The SRSM deals with this difficulty by approximating the output by a polynomial chaos expansion \cite{xiu2010numerical}. The below formulations follows \cite{huang_extended_2007}. The multidimensional Hermite polynomials of degree $p$ are used in the SRSM and defined as:
\begin{equation}
    \begin{aligned}
    &H_p(\alpha_{i_1},\alpha_{i_2},\ldots,\alpha_{i_p})\\
    &=(-1)^pe^{\frac{1}{2}\boldsymbol{\alpha}^T\boldsymbol{\alpha}}\frac{\partial^p}{\partial\alpha_{i_1},\partial\alpha_{i_2},\ldots,\partial\alpha_{i_p}}e^{-\frac{1}{2}\boldsymbol{\alpha}^T\boldsymbol{\alpha}}
    \end{aligned}
    \label{e:hermite}
\end{equation}
where $\boldsymbol{\alpha}=\{\alpha_{i_k}\}_{k=1}^p$ is a vector of standard normal random variables. The output of interest $z$ is estimated as follows:
\begin{equation}
    \begin{aligned}
    z = a_0 &+ \sum_{i_1=1}^n a_{i_1}H_1(\alpha_{i_1}) + \sum_{i_1=1}^n\sum_{i_2=1}^{i_1} a_{i_1i_2}H_2(\alpha_{i_1},\alpha_{i_2})\\
    &+\sum_{i_1=1}^n\sum_{i_2=1}^{i_1}\sum_{i_3=1}^{i_2} a_{i_1i_2i_3}H_3(\alpha_{i_1},\alpha_{i_2},\alpha_{i_3}) + \ldots
    \end{aligned}
    \label{e:srs}
\end{equation}
where $n$ is the number of standard normal random variables used in the expansion, and $a_0,a_{i_1},a_{i_1i_2},a_{i_1i_2i_3},\ldots$ are unknown coefficients. If $n=2$ and $p=3$, then the expansion~(\ref{e:srs}) will become:
\begin{equation}
    \begin{aligned}
    z(\alpha_{i_1},\alpha_{i_2})&=&a_0&+a_1\alpha_{i_1}+a_2\alpha_{i_2}+a_3(\alpha_{i_1}^2-1)+a_4(\alpha_{i_2}^2-1)\\
    && &+a_5\alpha_{i_1}\alpha_{i_2}+a_6(\alpha_{i_1}^3-3\alpha_{i_1})+a_7(\alpha_{i_2}^3-3\alpha_{i_2})\\
    && &+a_8(\alpha_{i_1}\alpha_{i_2}^2-\alpha_{i_1})+a_9(\alpha_{i_1}^2\alpha_{i_2}-\alpha_{i_2})\\
    &=&a_0&+\sum_{k=1}^9 a_kh_{i_k}
    \end{aligned}
    \label{e:n2p3}
\end{equation}
where $1,h_{i_1},h_{i_2},\ldots,h_{i_9}$ are Hermite polynomials. The ten unknown coefficients $a_0,a_1,\ldots,a_9$ are found by solving a system of linear equations using at least ten different realizations of $(\alpha_{i_1},\alpha_{i_2})$. Such realizations can be chosen at collocation points according to the SSG in Section~\ref{ssg}, or a much simpler heuristic rule in \cite{isukapalli_computationally_2004}, which is selected in this paper. For the approximation in~(\ref{e:n2p3}) the rule generates 17 collocation points to form a stochastic response surface, which was very close to the target output \cite{huang_extended_2007}.

\subsubsection{Solution Algorithm}
\label{algorithm}
The optimization algorithm (Fig.~\ref{fig:rr_sora}) to solve the RRBTO problem~(\ref{e:2}) is expounded below:
\begin{enumerate}
	\item Initialize the problem: size of the finite element mesh; initial values of design variables, SIMP and optimization parameters; KL expansion of random field; collocation points and weights for the SSG, and the SRSM; etc.
	\item Until convergence do:
	\setlist[itemize]{leftmargin=0mm}
    \begin{itemize}[label=\small$\bullet$]
    \item Solve $\text{\bf{K}}_i\text{\bf{u}}_i = \text{\bf{f}}_i$ and compute $\dfrac{\partial C_i(\boldsymbol{\rho})}{\partial \boldsymbol{\rho}}$ for $i=1,2,\ldots,n_l^{(d)}$.
	      	
  	\item Calculate mean, variance and their derivatives:
    \begin{equation}
    	\begin{aligned}
    		\mathds{E}\left[C\right]                                             & = \sum_{i=1}^{n_l^{(d)}} w_i C_i,                                                                                                                                                                                        \\
    		\sigma^2\left[C\right]                                               & = \sum_{i=1}^{n_l^{(d)}} w_i C_i^2 - \mathds{E}^2\left[C\right],                                                                                                                                                         \\
    		\frac{\partial \mathds{E}\left[C\right]}{\partial \boldsymbol{\rho}} & = \sum_{i=1}^{n_l^{(d)}} w_i \frac{\partial C_i(\boldsymbol{\rho})}{\partial \boldsymbol{\rho}},                                                                                                                         \\
    		\frac{\partial \sigma^2\left[C\right]}{\partial \boldsymbol{\rho}}   & = \sum_{i=1}^{n_l^{(d)}} 2C_i(\boldsymbol{\rho}) w_i \frac{\partial C_i(\boldsymbol{\rho})}{\partial \boldsymbol{\rho}} - 2\mathds{E}\left[C\right]\frac{\partial \mathds{E}\left[C\right]}{\partial \boldsymbol{\rho}}. 
    	\end{aligned}
    \end{equation}
    
  	\item Compute derivative of the robust objective:
    \begin{equation}
    	\kappa_1\frac{\partial \mathds{E}\left[C\right]}{\partial \boldsymbol{\rho}} + \kappa_2 \frac{1}{2\sqrt{\sigma^2[C]}} \frac{\partial \sigma^2\left[C\right]}{\partial \boldsymbol{\rho}}.
    \end{equation}
  	      
  	\item Deterministic topology optimization (DTO): the most probable point (MPP) $\boldsymbol{\xi}_i^{*(k-1)}$ found in the previous loop (or some initial values for the first loop) is used in place of random parameters $\boldsymbol{\xi}(\omega)$, making~(\ref{e:relia4}) a regular TO problem. The SIMP and the MMA method are employed to solve it.
  	\item Inverse reliability analysis (IRA): the optimum values of design variables from the DTO step and the collocation points given in Section~\ref{srsm} are used to construct stochastic response surfaces, which in turn are utilized in~(\ref{e:relia5}) to find the next MPP. Based on convergence conditions, the algorithm may stop, or a new loop is requested with updated Young's modulus using the new MPP.
  	\end{itemize}
\end{enumerate}

\section{Results}
\label{results}
In this section our proposed algorithm is run on two common benchmark problems (the cantilever and the L-shaped beam) with three target reliability levels, six weighting factors, and one parameter tuple of the uniform distribution in~(\ref{e:uniform}). The correctness and accuracy of our algorithm is then verified on the optimization results by Monte Carlo simulations. All quantities given below are dimensionless for simplicity.

To ensure the generality of our approach, we intentionally do not specify the limit state functions $g_i( \boldsymbol{\rho},y)$ in the previous sections, which is essential for the two numerical examples. The RRBTO problem now becomes:
\begin{equation}
    \begin{aligned}
        \min_{\boldsymbol{\rho}} && &\kappa_1 \mu\left[C( \boldsymbol{\rho},y )\right] + \kappa_2 \sigma\left[C( \boldsymbol{\rho},y )\right] \\
        \text{s.t.} && &\text{\bf{K}}(\boldsymbol{\rho},y) \text{\bf{u}}(\boldsymbol{\rho},y) = \text{\bf{f}}, \\
        && & \frac{V(\boldsymbol{\rho})}{V_0} = \gamma, \\
        && &P\left[u( \boldsymbol{\rho},y) - u_0 < 0\right] \leq P^0, \\
        && &0 < \rho_{\text{min}} \leq \boldsymbol{\rho} \leq 1,
        \label{e:examples}
    \end{aligned}
\end{equation}
where $u( \boldsymbol{\rho},y)$ and $u_0$ are the actual displacement and the minimum allowable displacement at a selected point, respectively. The above formulation is inspired by design of compliant mechanisms \cite{frecker_topological_1997}, in which both flexibility and stiffness are required. Flexibility allows the mechanisms to reach designed deformation, implied by the limit state function $g(\boldsymbol{\rho},y)=u(\boldsymbol{\rho},y) - u_0$, while maximizing stiffness, or minimizing compliance, helps them withstand loads. Furthermore, the weighting factors $\kappa_1$ and $\kappa_2$ need substantial attention. Our preliminary numerical results showed that the mean and standard deviation of the compliance are different by about two orders of magnitude, which may have critical impact on the solution. This was examined carefully in \cite{marler_weighted_2010} and normalization transforming them into the same scale has been recommended. $\kappa_1$ and $\kappa_2$ are identified according to \cite{asadpoure_robust_2011}:
\begin{equation}
    \kappa_1=\frac{\epsilon}{\mu^*}, \kappa_2=\frac{1-\epsilon}{\sigma^*},
    \label{e:k1k2}
\end{equation}
where $\epsilon \in [0,1]$, $\mu^*$ is the mean of the compliance when $(\epsilon,1-\epsilon)=(0,1)$, and $\sigma^*$ is the standard deviation of the compliance when $(\epsilon,1-\epsilon)=(1,0)$. $\mu^*$ and $\sigma^*$ have to be calculated on the same set of input parameters except $\epsilon$, under which they are the maximum values of the mean and standard deviation of the compliance resulting in $0 < \dfrac{\mu}{\mu^*},\dfrac{\sigma}{\sigma^*} \leq 1 $ in (\ref{e:examples}).

Coding the algorithm demands concrete values of every parameter, many of which are shared between the two examples. The finite element mesh in both examples is assembled from square, linear, plane stress elements, whose side length and thickness are unit dimension. Those elements are made of an isotropic, linear elastic material with Poisson's ratio $\nu=0.3$. The material Young's modulus is assumed to be a centered, mean-square Gaussian random field with known covariance function as in~(\ref{e:covariance}), expanded into a series of independent, standard normal random variables \cite{noauthor_brief_nodate}. The correlation lengths are chosen as $l_1=l_2=0.6$, and only the first two eigenvalues and eigenfunctions are picked for the truncated KL expansion. A wide range of optimization algorithms for reliability analysis is available in the CODES toolbox \cite{codes_toolbox}, including the Hybrid Mean Value method selected for the examples due to its efficiency \cite{youn_hybrid_2003}. Other parameters of the optimization problem (\ref{e:examples}) are the target reliability levels $\beta=\{1.0,2.0,3.0\}$, weighting factors $\epsilon=\{1,0.9,0.8,0.5,0.2,0\}$, and the material property limits $(a,b) = (1,1.5)$. Due to iterative nature of MMA and SORA, several convergence criteria are enforced. The MMA optimizer stops when the maximum difference of design variables of two consecutive iterations is smaller than a prescribed value ($d_{\text{max}}^{\text{MMA}} \leq 0.001$), or the number of iterations is more than $n_{\text{MMA}}=200$. The similar conditions are applied in the SORA loops, but for the MPPs ($d_{\text{max}}^{\text{MPP}} \leq 0.001$) and a different maximum allowable number of iterations ($n_{\text{SORA}}=20$). The level of sparse grid approximation is 4 using Gauss-Hermite quadrature as the base one-dimensional rule. In SIMP, the minimum length scale $r_{\text{min}}=1.5$ and penalization factor $p=3$ are used. Each example is then validated by 50000 Monte Carlo simulations (MCS), which compute failure probabilities of the limit state function, and the mean and standard deviation of the compliance. The examples are implemented in Matlab using Latin hypercube sampling for MCS with seed 0. Those probabilities are compared with values calculated from the three reliability levels, while statistical moments of the compliance prove robustness of the optimization results against uncertainty.

\subsection{The Cantilever Beam}
\label{cantilever}
\begin{figure}[H]
    \centering
    \includegraphics[width=0.35\textwidth]{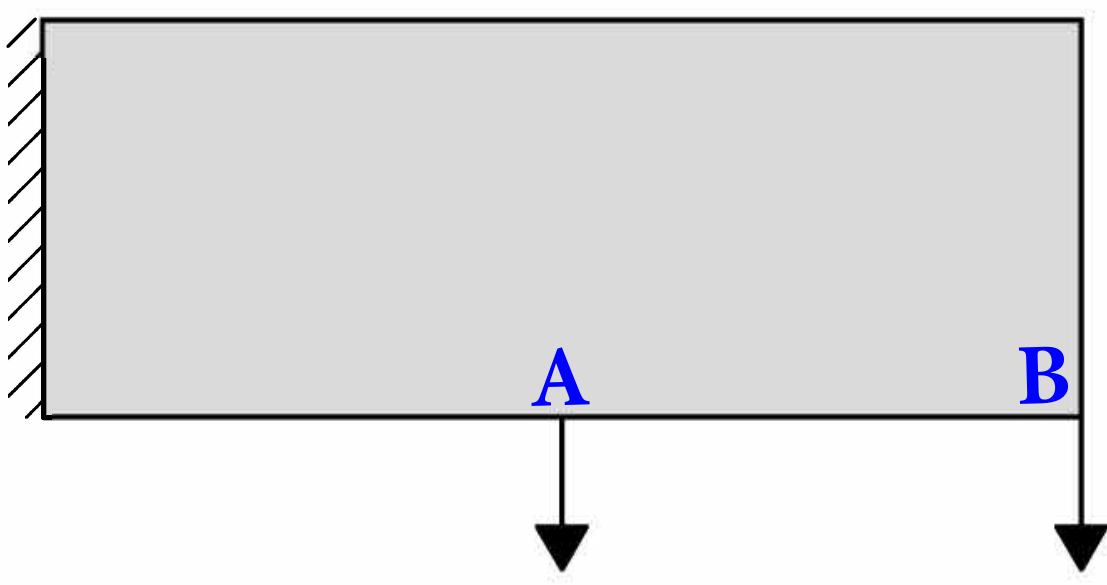}
    \caption{The cantilever beam}
    \label{fig:C}
\end{figure}

Fig.~\ref{fig:C} displays the cantilever beam, a two-dimensional domain used in this example. The beam is fixed on its left side, meshed into 60 $\times$ 20 elements, and subject to two unit vertical loads, which are placed on its bottom edge at equal distances (points A and B). The optimization problem is set up as in~(\ref{e:examples}), in which the minimum allowable displacement at the load application point B is given as $u_0=220$. Then our proposed algorithm is tested on the example with different values of target reliability, material parameters and weighting factors. The 18 optimized designs are presented in Tables \ref{tbl:RRBTO_C_1} and \ref{tbl:RRBTO_C_2}. The MCS and SRSM utilize those designs to calculate the mean and standard deviation of vertical displacement of point B ($\mu_B$ and $\sigma_B$ in Table \ref{tbl:C_Data}), as well as the probabilities of failure event $P_f=P\left[g(\boldsymbol{\rho},y) < 0\right]$ to show the reliability levels achieved by the designs. Moreover, the mean and standard deviation of the compliance are also obtained from the MCS ($\mu[C]$ and $\sigma[C]$ in Table \ref{tbl:C_Data}) to examine how they vary with respect to the weighting factors $\epsilon$. All of theses values are gathered in Table \ref{tbl:C_Data}, whose second column shows expected failure probabilities $P_f=\Phi(-\beta)$. Some conclusions from this example results can be found in Section~\ref{disu}.

\begin{sidewaystable*}[ph!]
    \centering
    \caption{The cantilever beam: RRBTO results}
    \vspace{5px}
    \begin{tabular}{lcccc}
        \toprule
        $(\epsilon,1-\epsilon)$ & $(1,0)$ & $(0.9,0.1)$ & $(0.8,0.2)$ \\
        \toprule
        $\beta=1$
        & \includegraphics[scale=0.5]{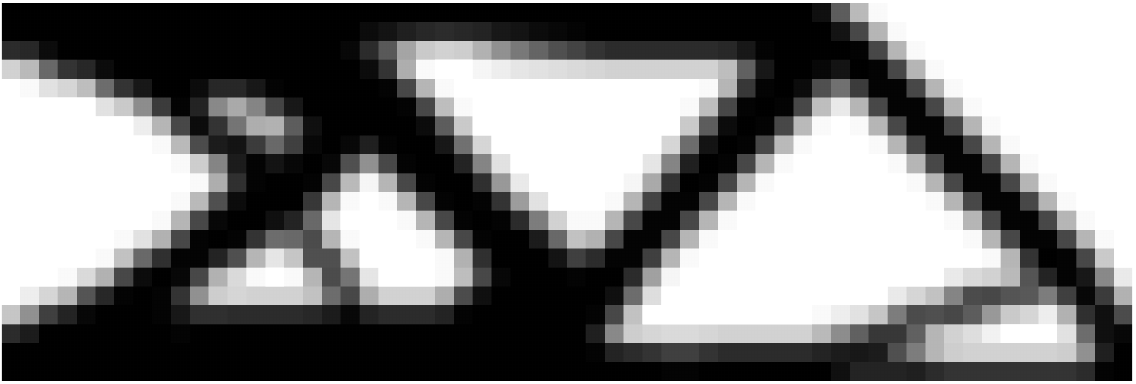}
        & \includegraphics[scale=0.5]{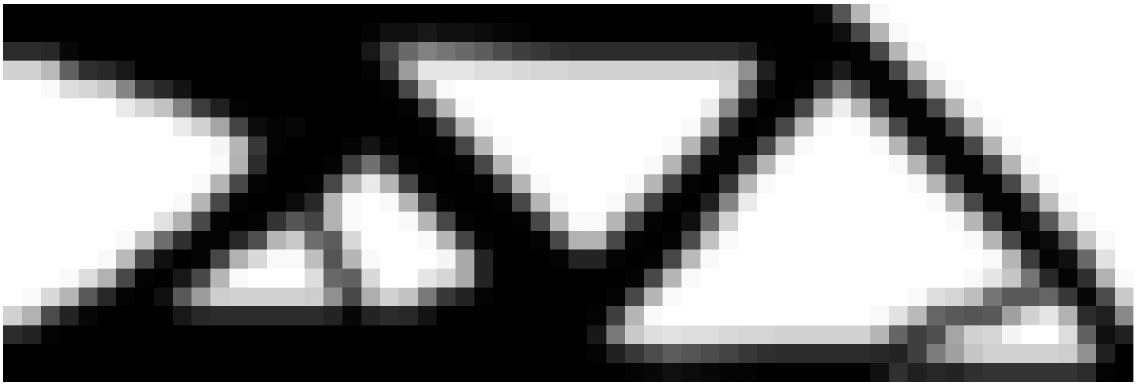}
        & \includegraphics[scale=0.5]{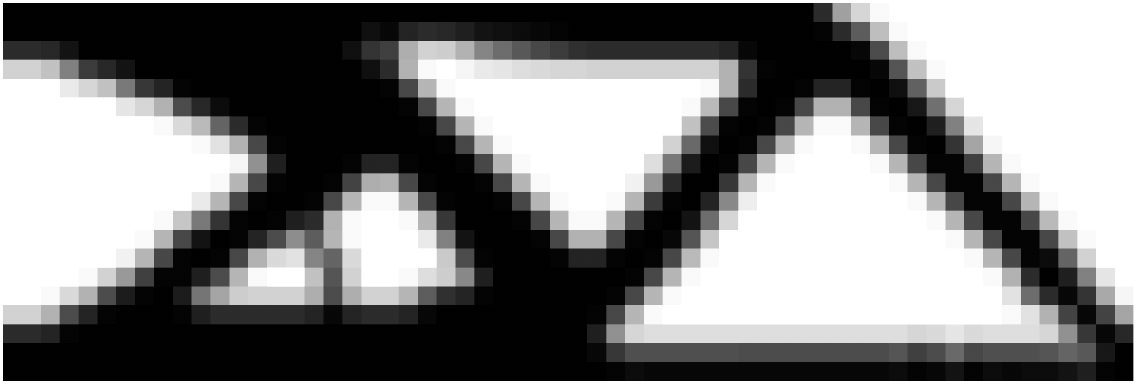} \\
        \midrule
        $\beta=2$ 
        & \includegraphics[scale=0.5]{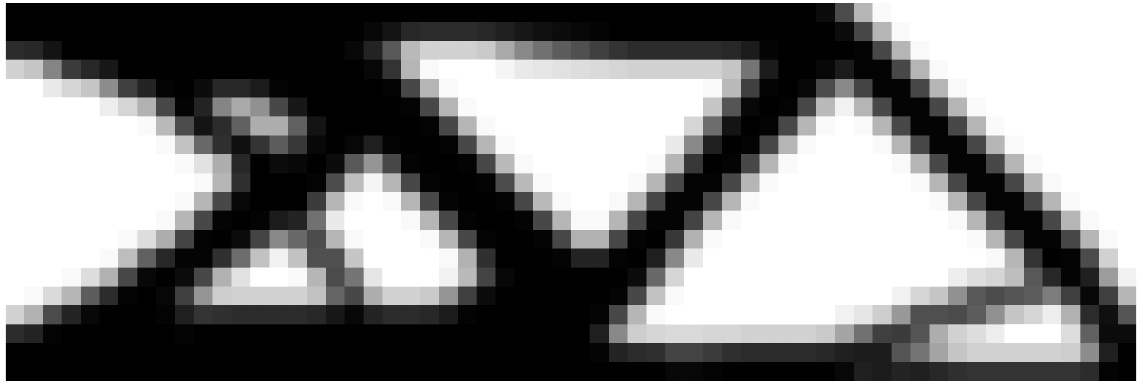}
        & \includegraphics[scale=0.5]{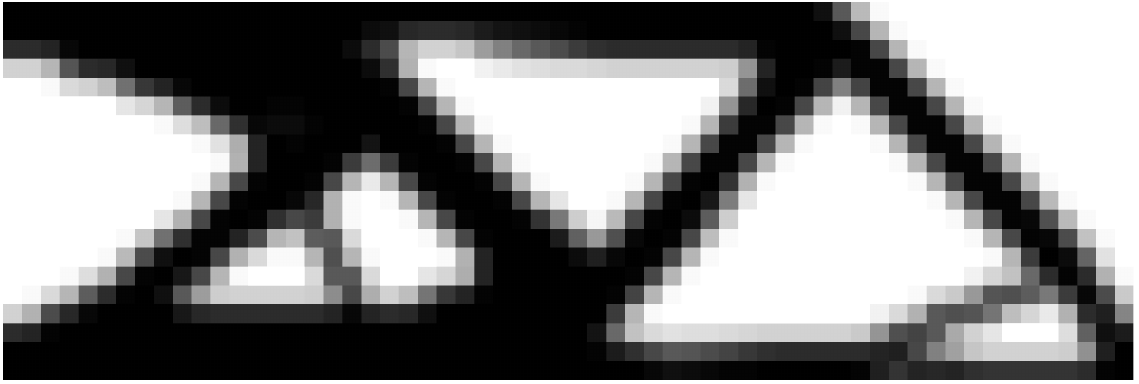} 
        & \includegraphics[scale=0.5]{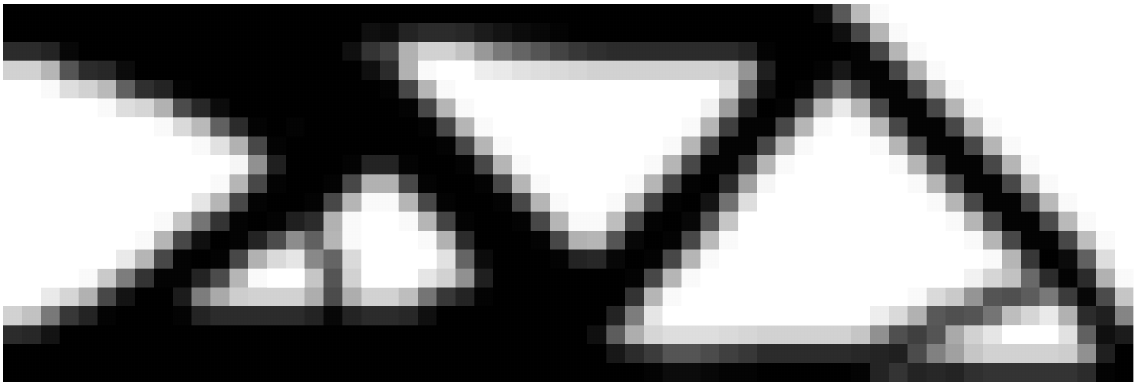} \\
        \midrule
        $\beta=3$ 
        & \includegraphics[scale=0.5]{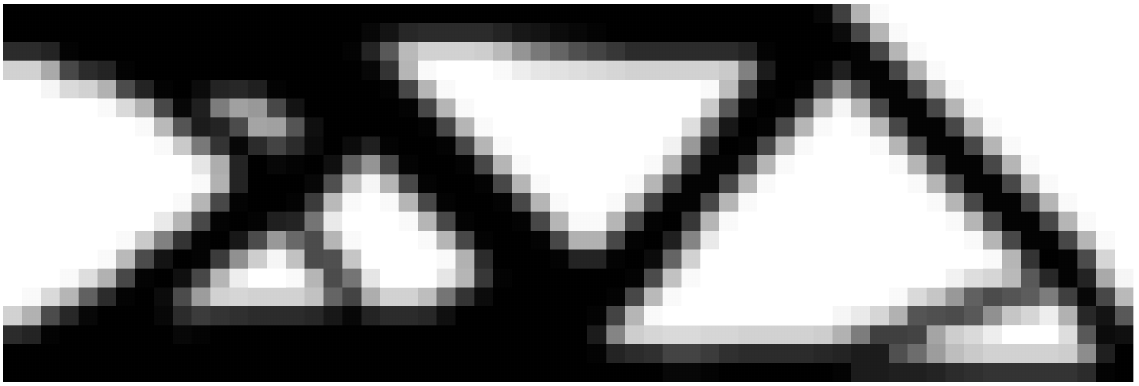}
        & \includegraphics[scale=0.5]{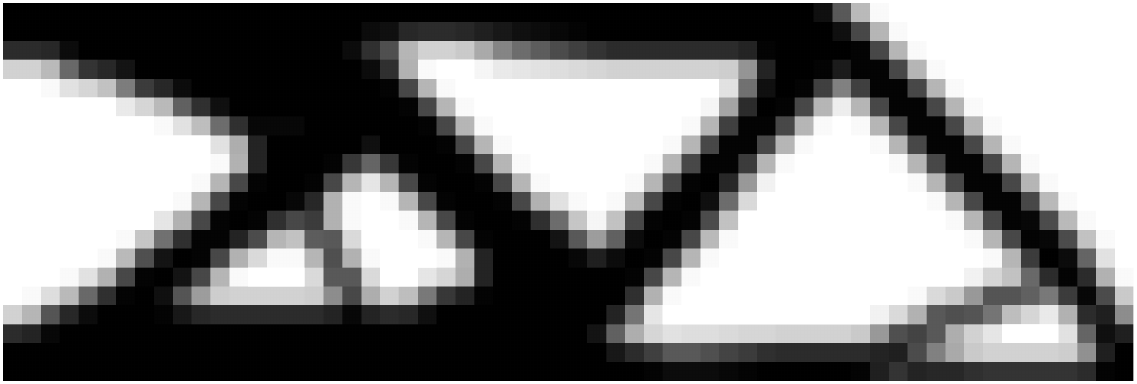}
        & \includegraphics[scale=0.5]{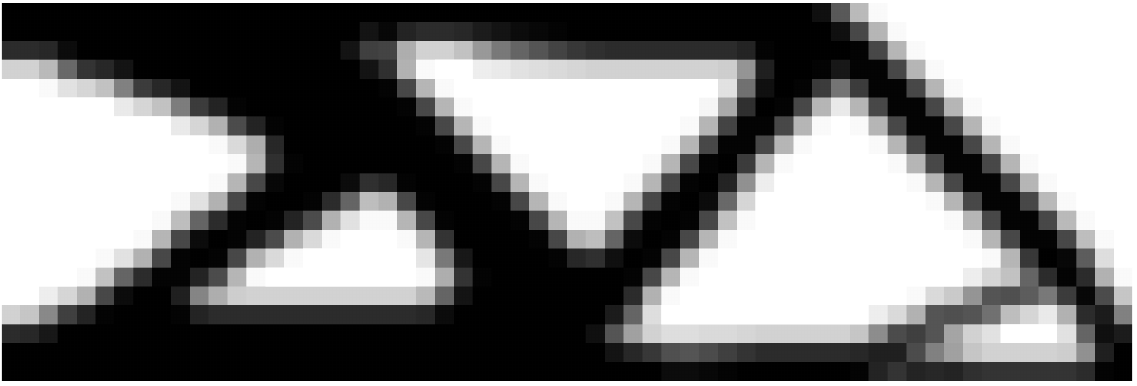} \\
        \bottomrule
        \label{tbl:RRBTO_C_1}
    \end{tabular}
    \caption{The cantilever beam: RRBTO results}
    \vspace{5px}
    \begin{tabular}{lcccc}
        \toprule
        $(\epsilon,1-\epsilon)$ & $(0.5,0.5)$ & $(0.2,0.8)$ & $(0,1)$ \\
        \toprule
        $\beta=1$
        & \includegraphics[scale=0.5]{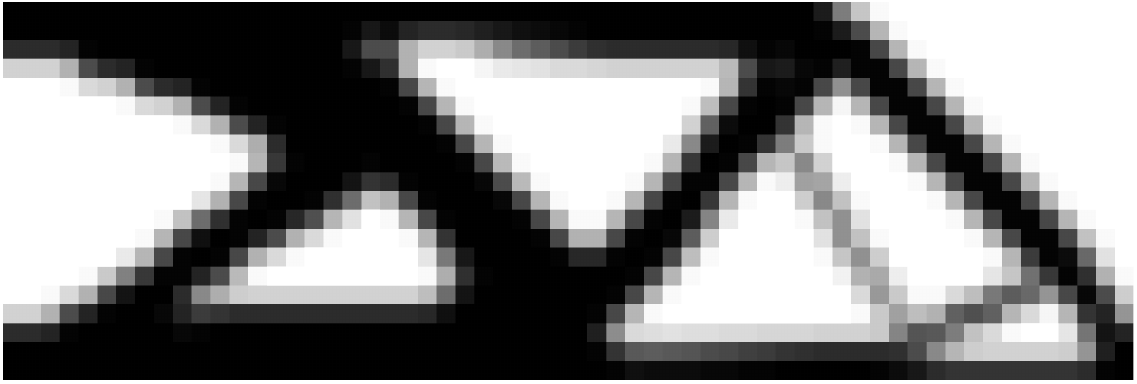}
        & \includegraphics[scale=0.5]{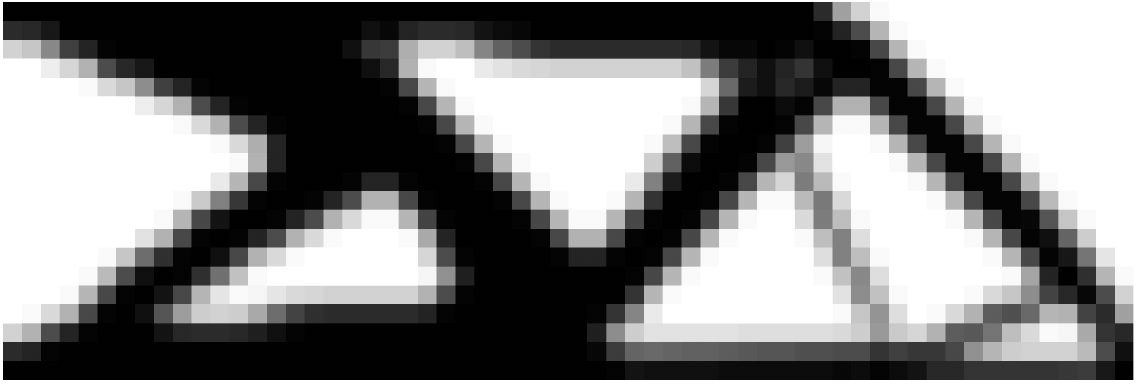}
        & \includegraphics[scale=0.5]{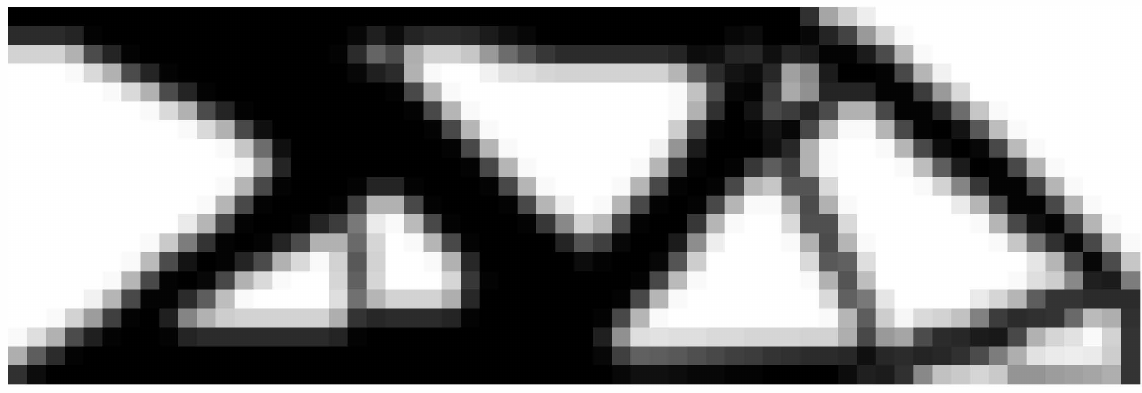} \\
        \midrule
        $\beta=2$ 
        & \includegraphics[scale=0.5]{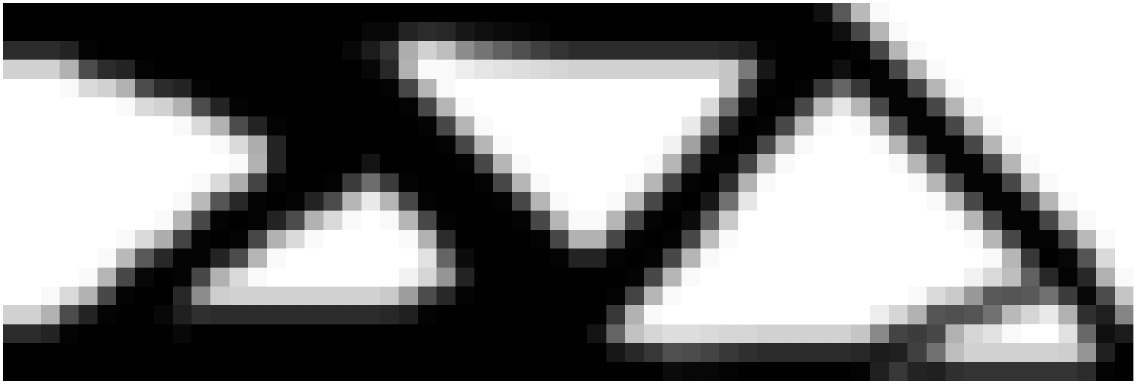}
        & \includegraphics[scale=0.5]{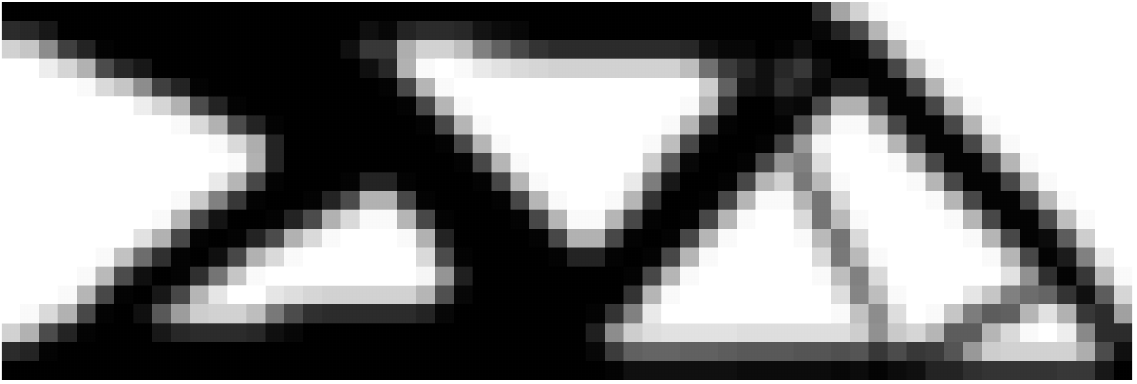} 
        & \includegraphics[scale=0.5]{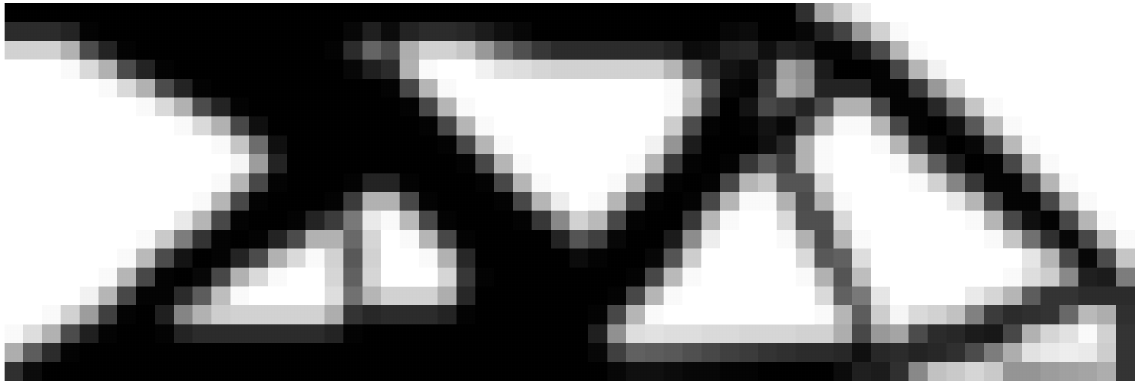} \\
        \midrule
        $\beta=3$ 
        & \includegraphics[scale=0.5]{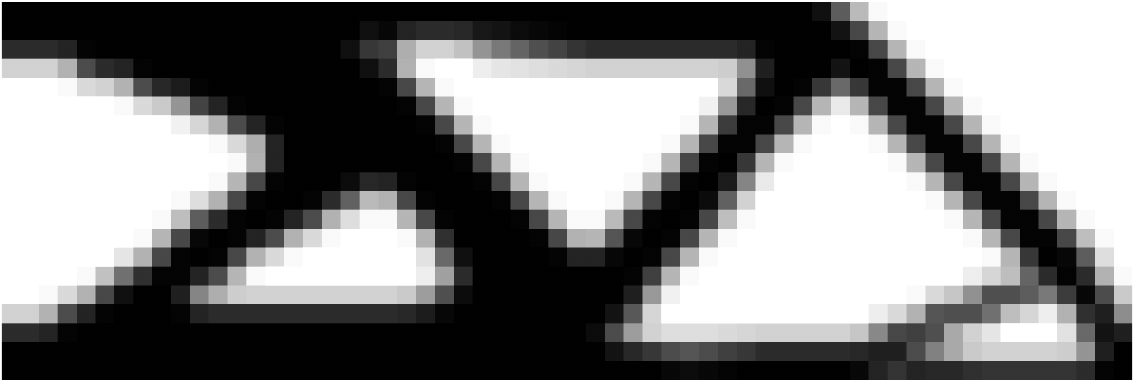}
        & \includegraphics[scale=0.5]{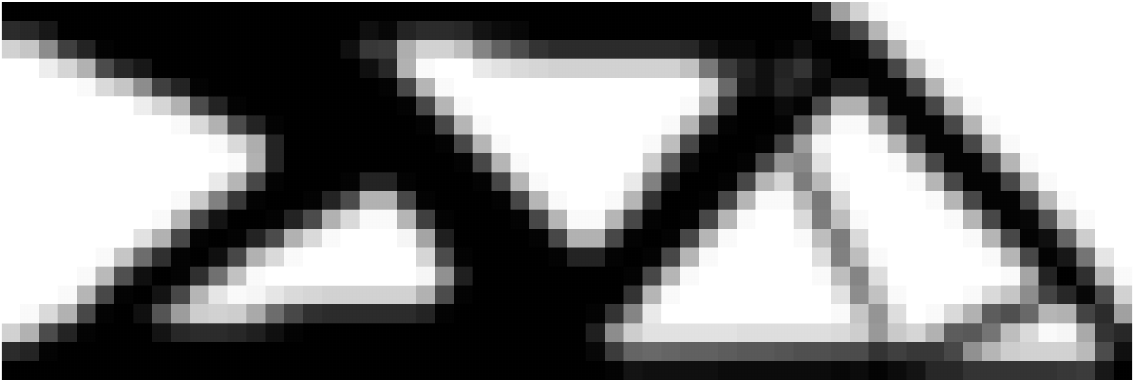}
        & \includegraphics[scale=0.5]{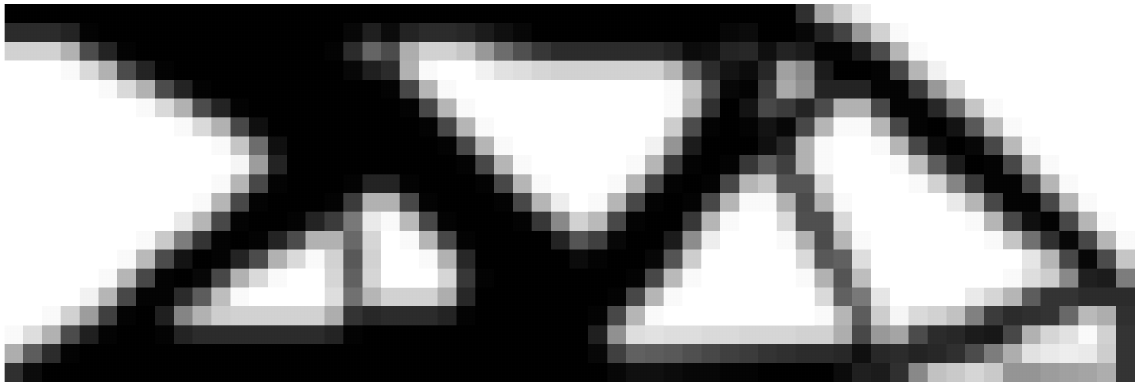} \\
        \bottomrule
        \label{tbl:RRBTO_C_2}
    \end{tabular}
\end{sidewaystable*}

\begin{table*}
    \centering
    \caption{The cantilever beam: Numerical results}
    \vspace{3px}
    \begin{tabular}{llllllllll}
        \toprule
        & & & & & \multicolumn{3}{c}{MCS} & \multicolumn{2}{c}{SRSM} \\
        \cmidrule(lrr){6-8}
        \cmidrule(lr){9-10}
        $\beta$ & Expected $P_f$ & $(\epsilon,1-\epsilon)$ & $\mu[C]$ & $\sigma[C]$ & $\mu_B$ & $\sigma_B$ & $P_f$ & $\mu_B$ & $\sigma_B$ \\
        \cmidrule{1-10}
        \multirow{6}{*}[-1em]{1} & \multirow{6}{*}[-1em]{0.15865} & $(1,0)$ & 162.9505&0.9263&-220.6120&0.6071&0.15674&-220.6120&0.6071 \\
        \cmidrule{3-10}
        & & $(0.9,0.1)$ & 162.9992&0.9188&-220.6070&0.6017&0.15642&-220.6070&0.6017 \\
        \cmidrule{3-10}
        & & $(0.8,0.2)$ & 163.2204&0.8953&-220.5980&0.5853&0.15326&-220.5980&0.5853 \\
        \cmidrule{3-10}
        & & $(0.5,0.5)$ & 166.5160&0.9032&-225.0110&0.5906&0.00000&-225.0110&0.5906 \\
        \cmidrule{3-10}
        & & $(0.2,0.8)$ & 177.7632&0.8588&-240.6990&0.5572&0.00000&-240.6990&0.5572 \\
        \cmidrule{3-10}
        & & $(0,1)$ & 206.2317&0.8559&-278.9100&0.5527&0.00000&-278.9100&0.5527 \\
        \cmidrule{1-10}
        \multirow{6}{*}[-1em]{2} & \multirow{6}{*}[-1em]{0.02275} & $(1,0)$ & 163.4230&0.9150&-221.1960&0.5995&0.02252&-221.1960&0.5995 \\
        \cmidrule{3-10}
        & & $(0.9,0.1)$ & 163.6040&0.9106&-221.1880&0.5961&0.02250&-221.1880&0.5961 \\
        \cmidrule{3-10}
        & & $(0.8,0.2)$ & 163.7847&0.9118&-221.1890&0.5967&0.02266&-221.1890&0.5967 \\
        \cmidrule{3-10}
        & & $(0.5,0.5)$ & 164.7008&0.8953&-222.9560&0.5866&0.00000&-222.9560&0.5866 \\
        \cmidrule{3-10}
        & & $(0.2,0.8)$ & 178.0176&0.8586&-241.0670&0.5570&0.00000&-241.0670&0.5570 \\
        \cmidrule{3-10}
        & & $(0,1)$ & 206.2320&0.8559&-278.9110&0.5527&0.00000&-278.9110&0.5527 \\
        \cmidrule{1-10}
        \multirow{6}{*}[-1em]{3} & \multirow{6}{*}[-1em]{0.001349} & $(1,0)$ & 163.7490&0.9094&-221.7690&0.5962&0.001340&-221.7690&0.5962 \\
        \cmidrule{3-10}
        & & $(0.9,0.1)$ & 164.0019&0.9078&-221.7630&0.5944&0.001340&-221.7630&0.5944 \\
        \cmidrule{3-10}
        & & $(0.8,0.2)$ & 164.0458&0.9015&-221.7520&0.5907&0.001320&-221.7520&0.5907 \\
        \cmidrule{3-10}
        & & $(0.5,0.5)$ & 165.0556&0.8938&-223.2560&0.5855&0.000000&-223.2560&0.5855 \\
        \cmidrule{3-10}
        & & $(0.2,0.8)$ & 177.8977&0.8585&-240.8960&0.5570&0.000000&-240.8960&0.5570 \\
        \cmidrule{3-10}
        & & $(0,1)$ & 206.2320&0.8559&-278.9110&0.5527&0.000000&-278.9110&0.5527 \\
        \bottomrule
        \label{tbl:C_Data}
    \end{tabular}
\end{table*}

\subsection{The L-shaped Beam}
\label{L-shaped}
\begin{figure}[H]
    \centering
    \includegraphics[width=0.27\textwidth]{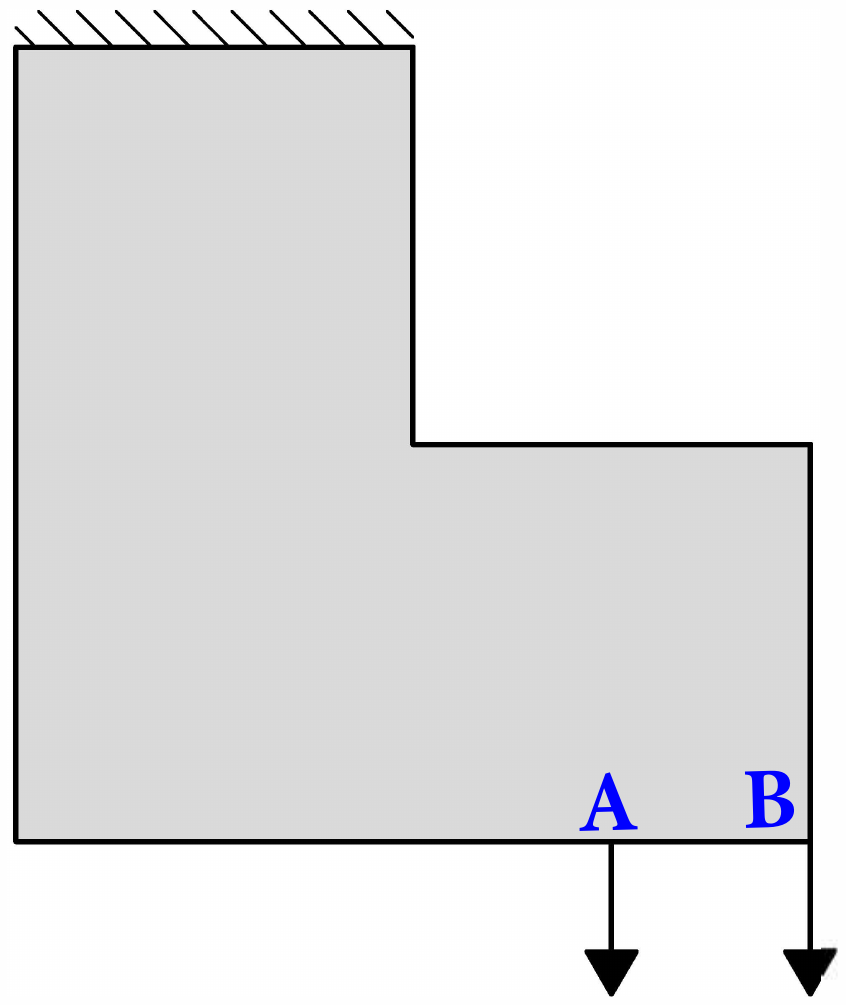}
    \caption{The L-shaped beam}
    \label{fig:L}
\end{figure}
The second numerical example to illustrate our proposed algorithm is the L-shaped beam as shown in Fig.~\ref{fig:L}. The beam is fixed on its topmost edge and subject to two unit vertical loads on its bottom edge$-$one at the right endpoint B and the other at point A located one quarter of the bottom edge length from point B. As in the previous example, the weighted sum of the two statistical moments of the compliance is minimized, while a probabilistic constraint is imposed on the vertical displacement of one load application point (point B in Fig.~\ref{fig:L}) to design for flexibility. The minimum allowable vertical displacement of point B is chosen as $u_0=130$. The design domain is discretized using a $60\times60$ mesh of finite elements, and then one quarter of the mesh is removed to make the domain L-shaped by forcing the element densities in this region to be $0.001$ before proceeding to the next operation. The larger mesh used in this examples results in much longer running time based on our experiments on the same computer. The results also exhibit the same trends as in the previous example. Therefore, instead of running the complete set of 18 combinations, only 12 cases, which are the combinations of $\beta=\{1,3\}$, six weighting factors, and $(a,b)=(1,1.5)$, are considered. Table~\ref{tbl:L_Data} shows the probabilities of displacement constraint violation $P_f$ at point B, the statistical moments of that point's vertical displacement ($\mu_B$ and $\sigma_B$) calculated from both the MCS and SRSM, and the mean and standard deviation of the compliance ($\mu[C]$ and $\sigma[C]$). The L-shaped optimized designs are analyzed for insights in the next section.

\begin{table*}[ph!]
    \centering
    \caption{The L-shaped beam: RRBTO results}
    \vspace{3px}
    \begin{tabular}{lcccc}
        \toprule
        $(\epsilon,1-\epsilon)$ & $(1,0)$ & $(0.9,0.1)$ & $(0.8,0.2)$ \\
        \toprule
        $\beta=1$
        & \includegraphics[scale=0.5]{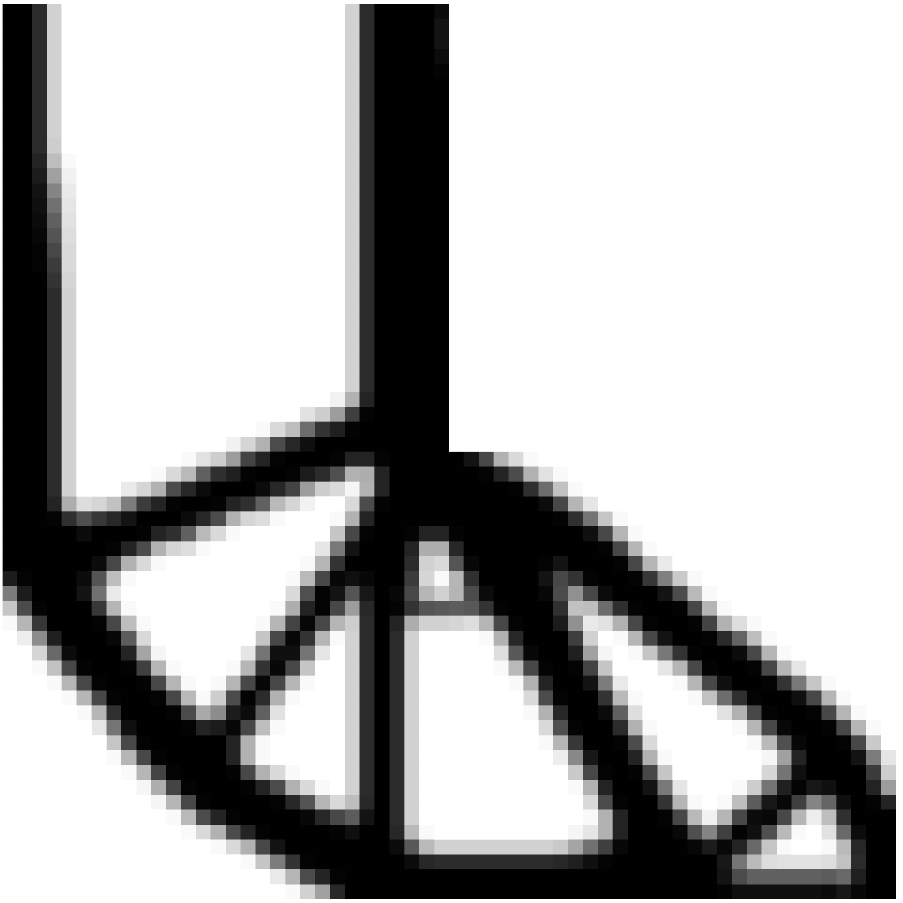}
        & \includegraphics[scale=0.5]{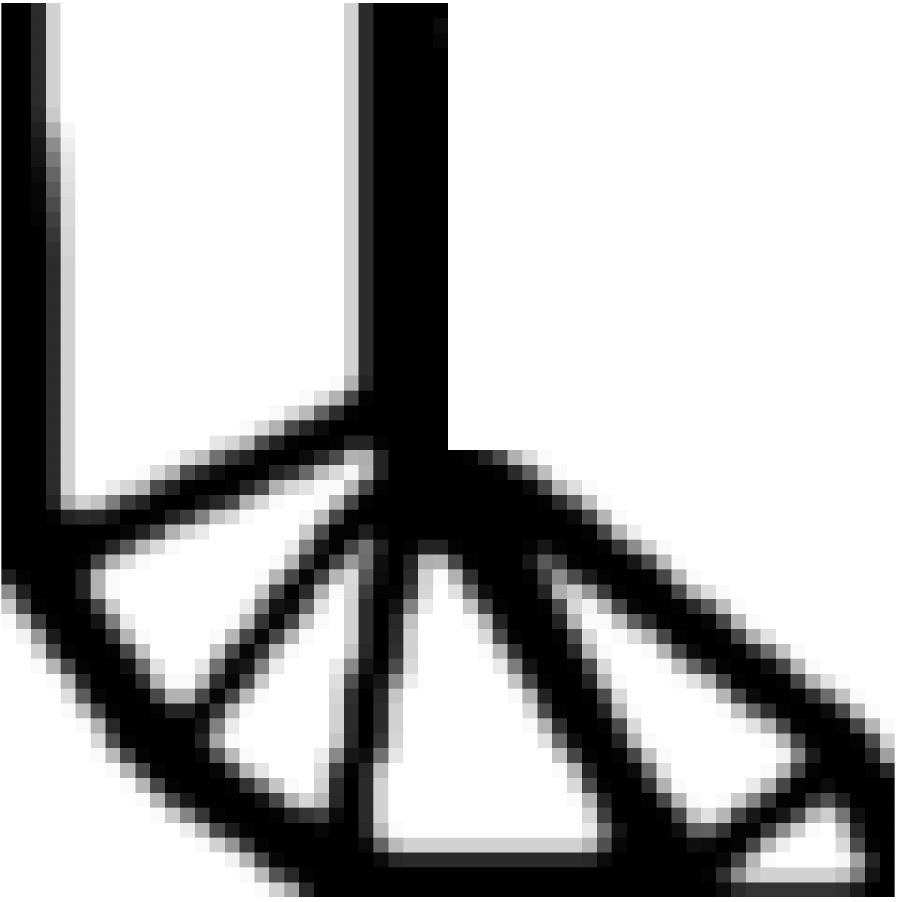}
        & \includegraphics[scale=0.5]{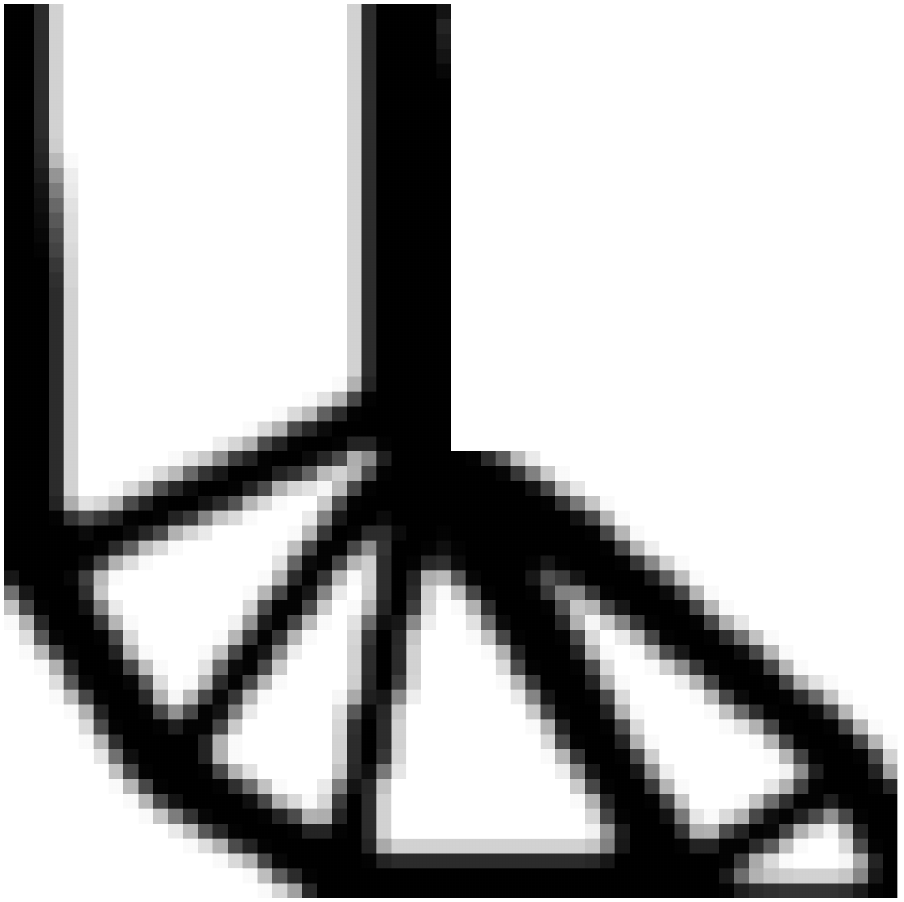} \\
        \midrule
        $\beta=3$ 
        & \includegraphics[scale=0.5]{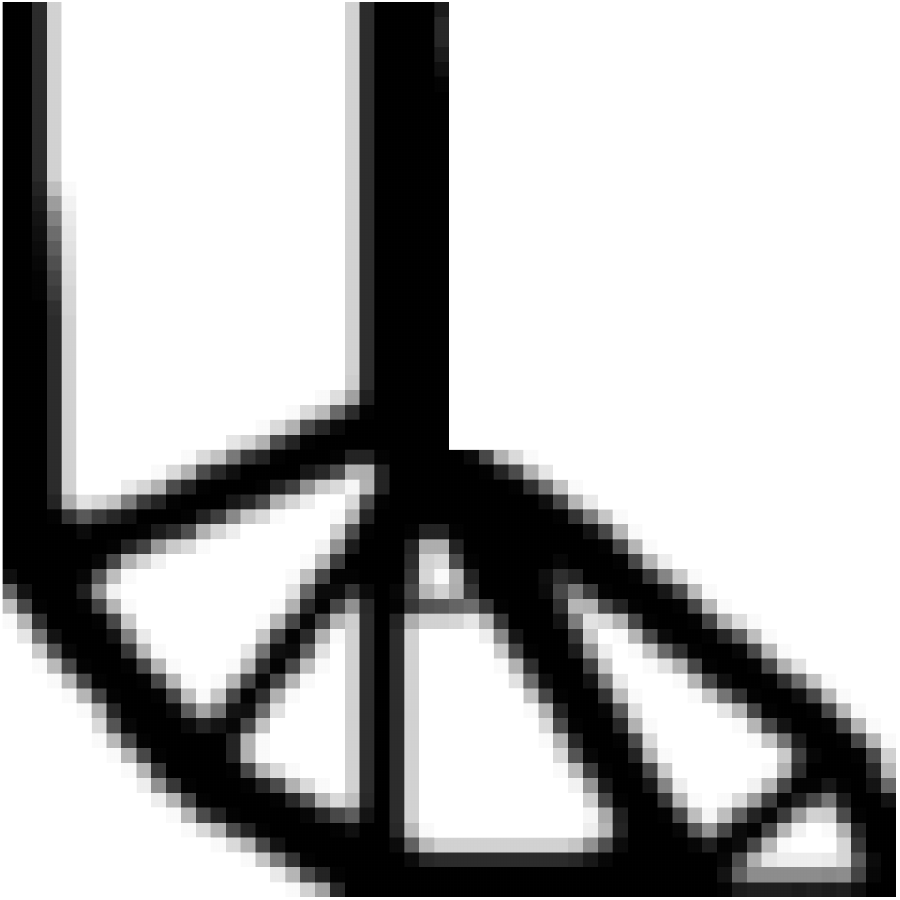}
        & \includegraphics[scale=0.5]{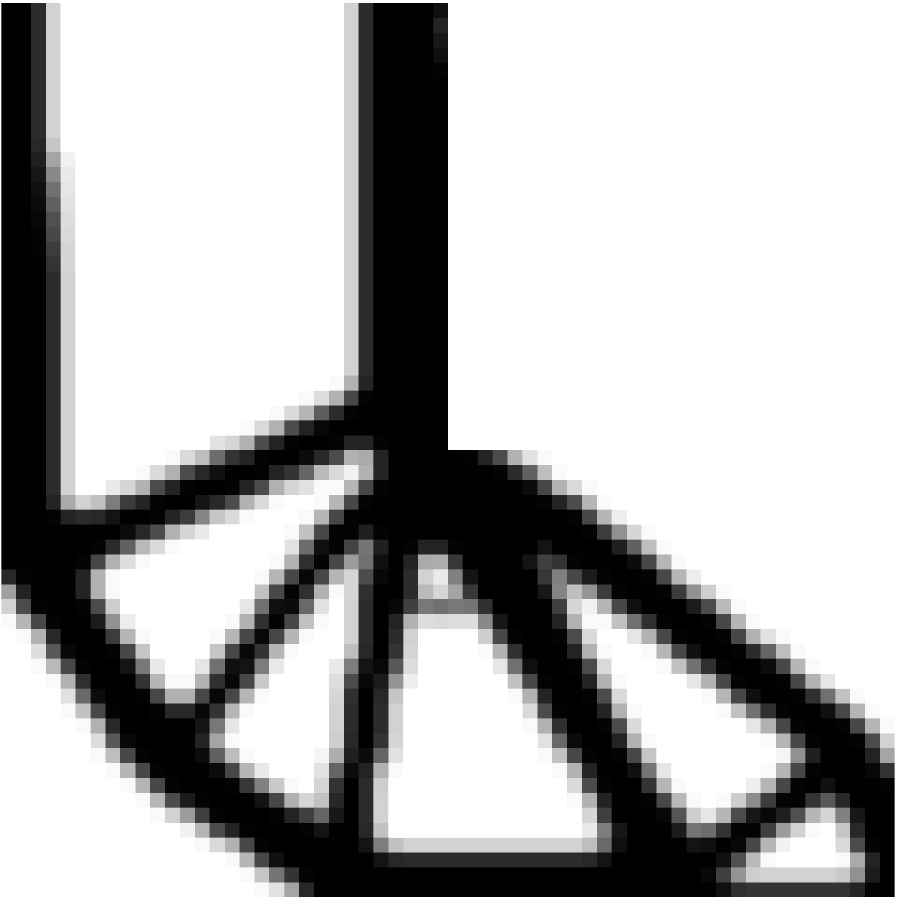}
        & \includegraphics[scale=0.5]{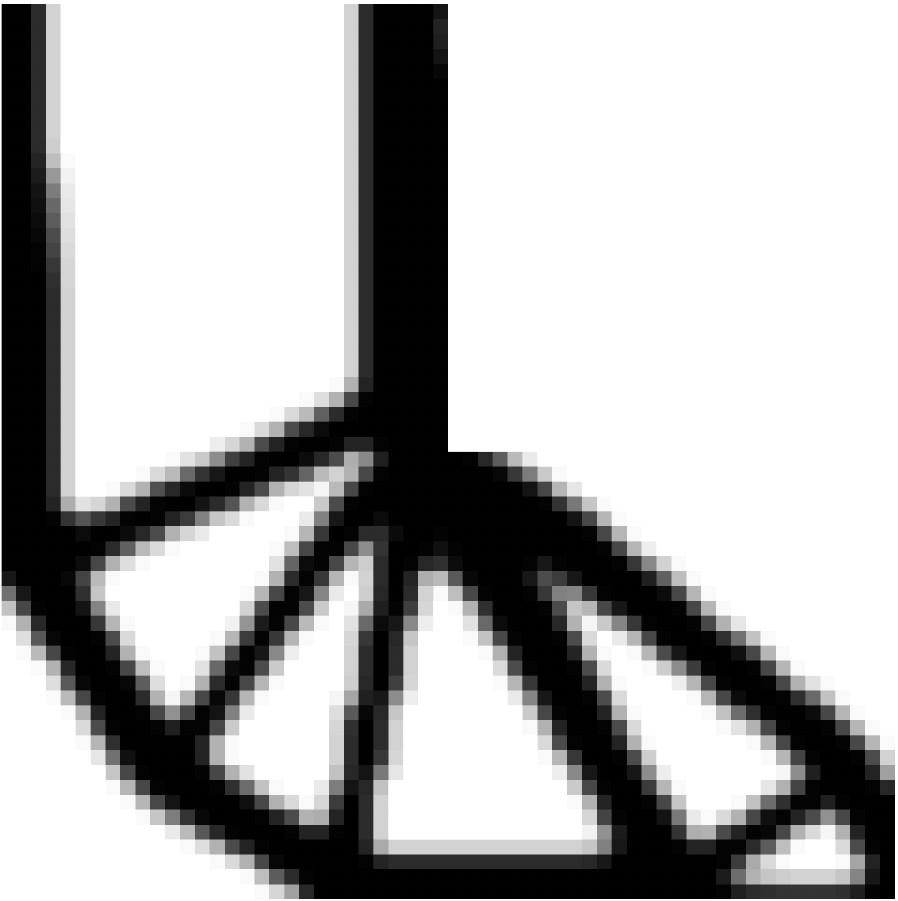} \\
        \bottomrule
        \label{tbl:RRBTO_L_1}
    \end{tabular}
    \caption{The L-shaped beam: RRBTO results}
    \vspace{3px}
    \begin{tabular}{lcccc}
        \toprule
        $(\epsilon,1-\epsilon)$ & $(0.5,0.5)$ & $(0.2,0.8)$ & $(0,1)$ \\
        \toprule
        $\beta=1$
        & \includegraphics[scale=0.5]{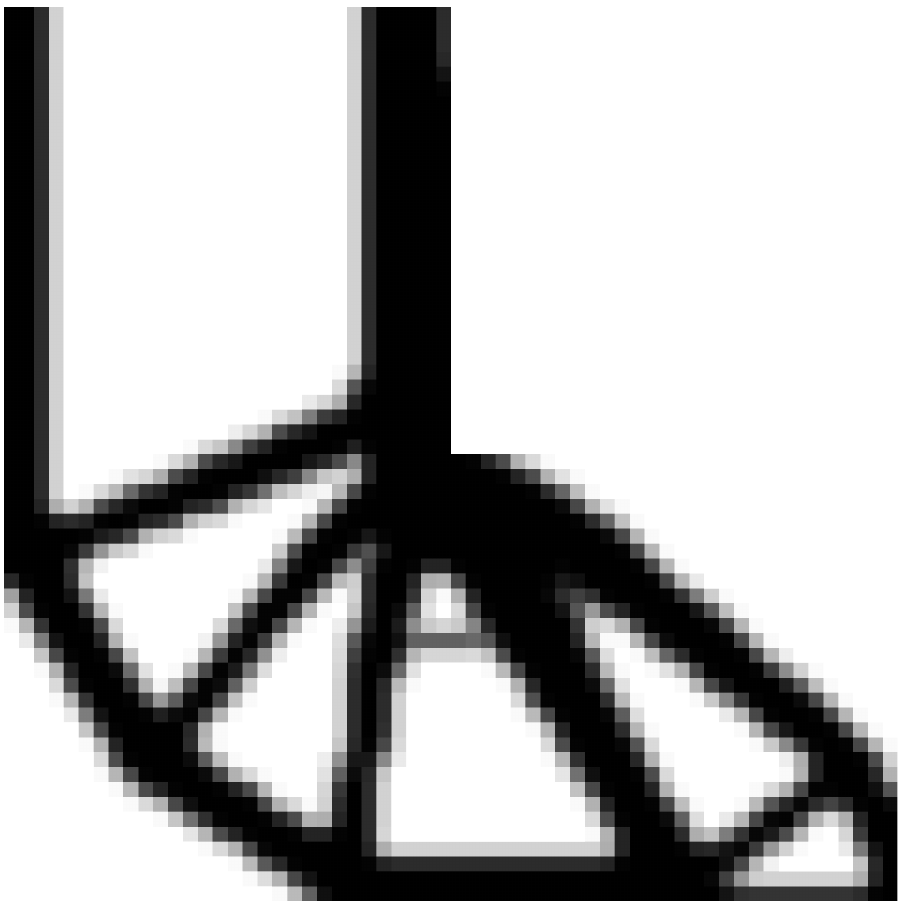}
        & \includegraphics[scale=0.5]{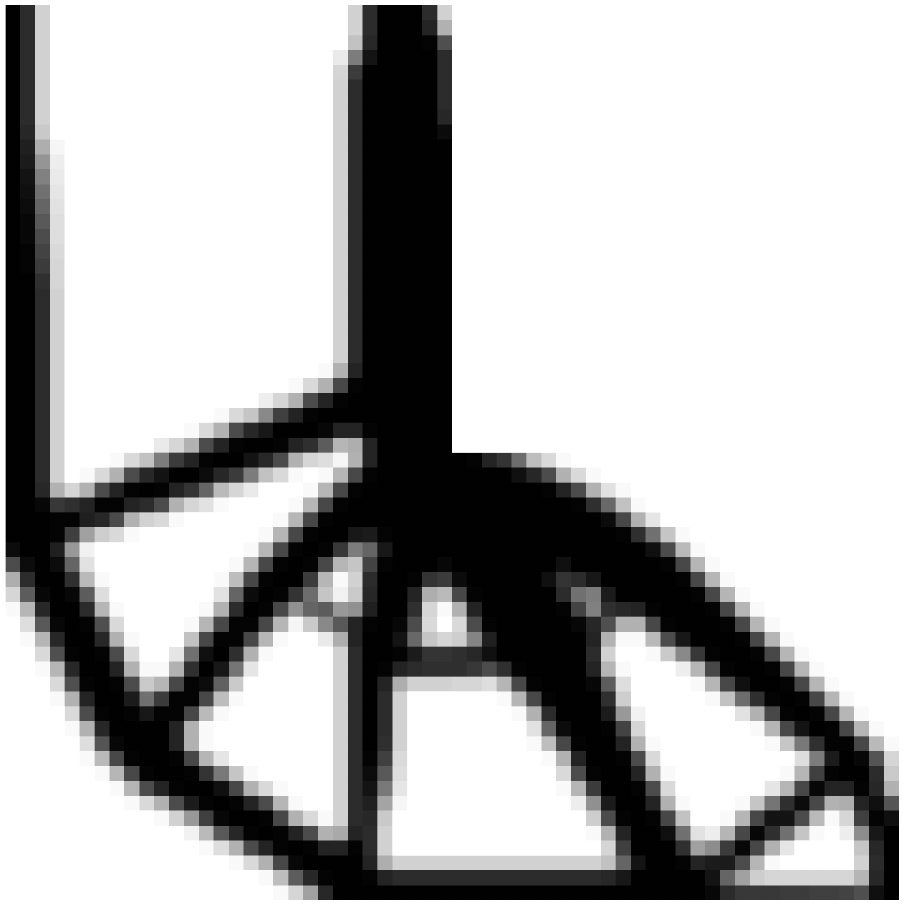}
        & \includegraphics[scale=0.5]{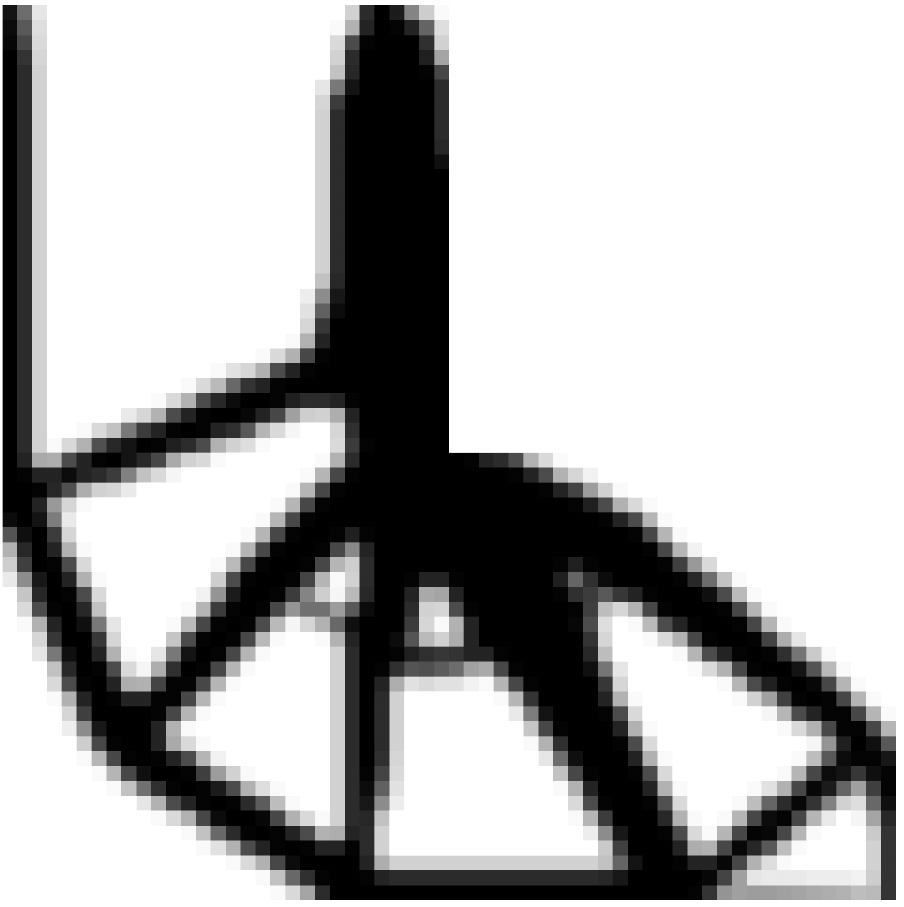} \\
        \midrule
        $\beta=3$ 
        & \includegraphics[scale=0.5]{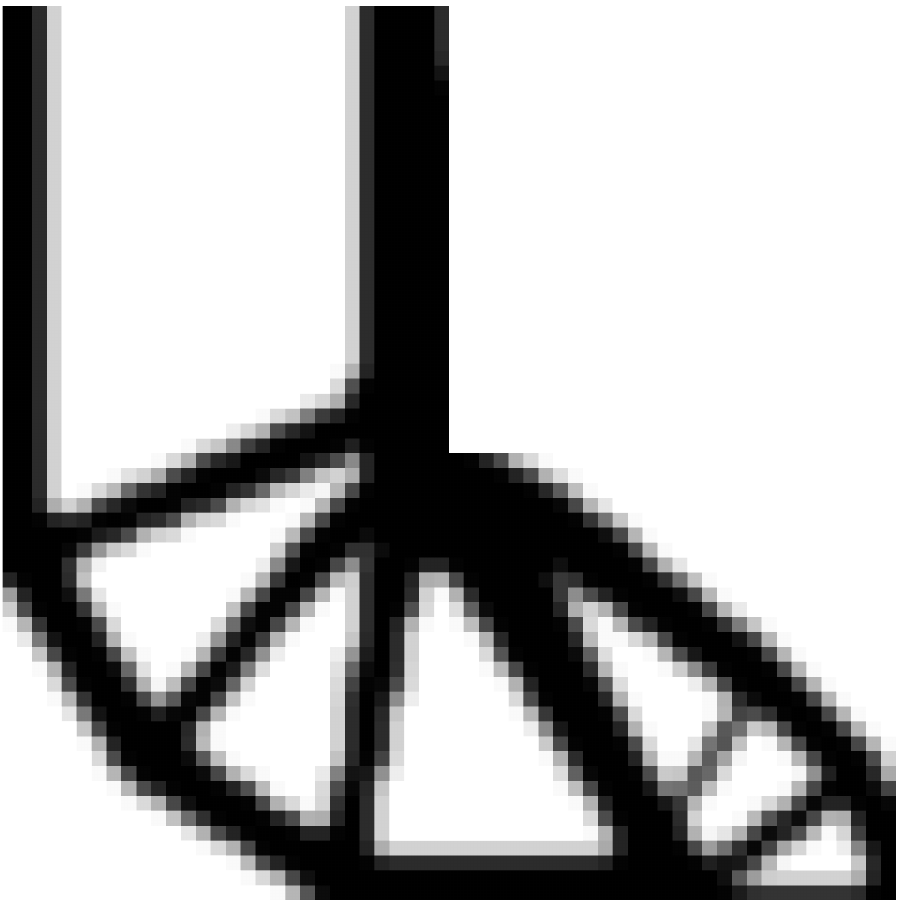}
        & \includegraphics[scale=0.5]{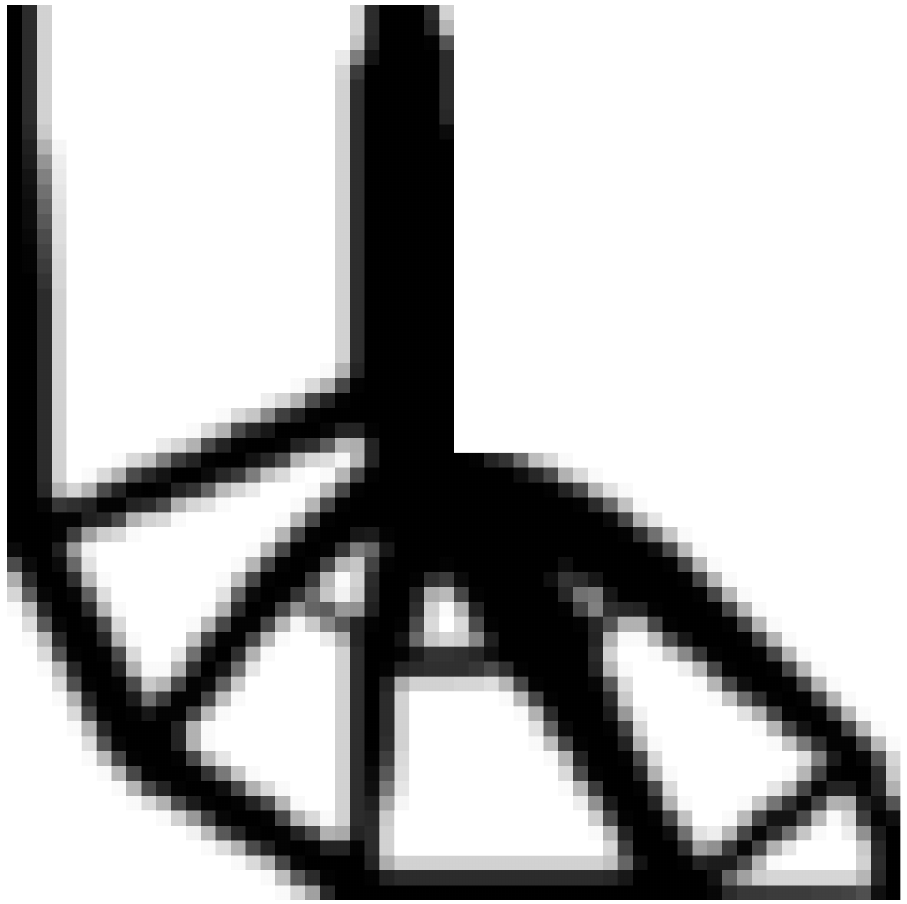}
        & \includegraphics[scale=0.5]{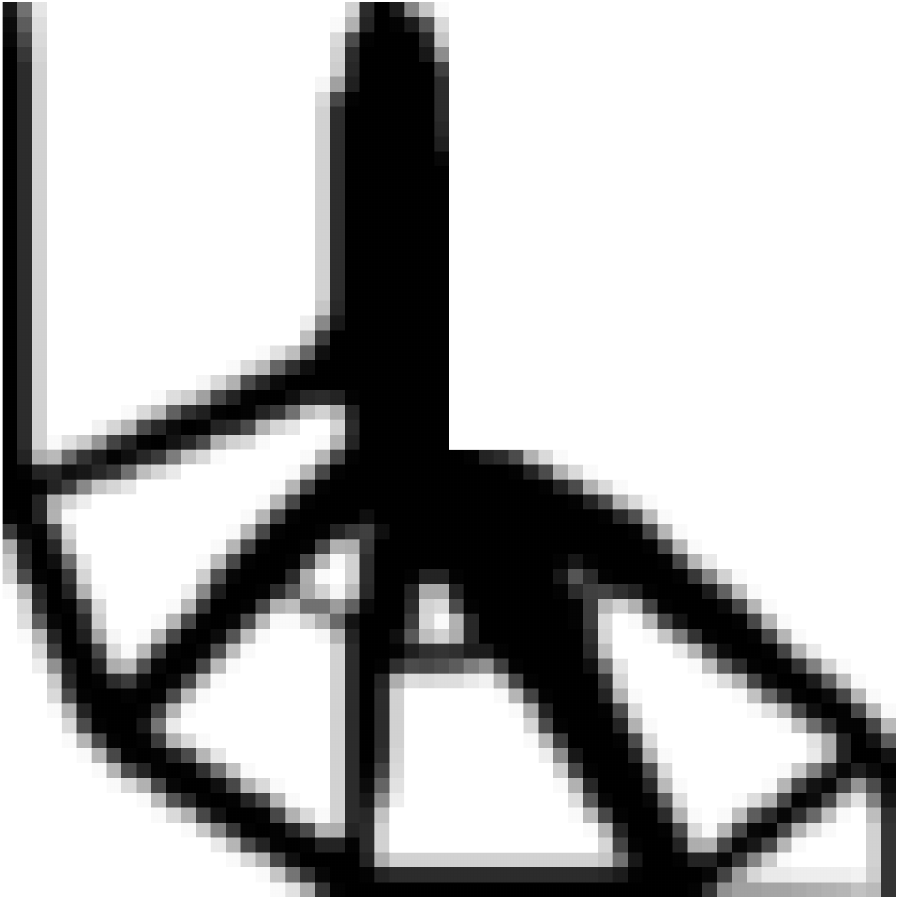} \\
        \bottomrule
        \label{tbl:RRBTO_L_2}
    \end{tabular}
\end{table*}

\begin{table*}
    \centering
    \caption{The L-shaped beam: Numerical results}
    \vspace{3px}
    \begin{tabular}{llllllllll}
        \toprule
        & & & & & \multicolumn{3}{c}{MCS} & \multicolumn{2}{c}{SRSM} \\
        \cmidrule(lrr){6-8}
        \cmidrule(lr){9-10}
        $\beta$ & Expected $P_f$ & $(\epsilon,1-\epsilon)$ & $\mu[C]$ & $\sigma[C]$ & $\mu_B$ & $\sigma_B$ & $P_f$ & $\mu_B$ & $\sigma_B$ \\
        \cmidrule{1-10}
        \multirow{6}{*}[-1em]{1} & \multirow{6}{*}[-1em]{0.15865} & $(1,0)$ & 96.4893&0.5352&-130.3580&0.3577&0.15808&-130.3580&0.3577 \\
        \cmidrule{3-10}
        & & $(0.9,0.1)$ & 96.3536&0.5321&-130.3560&0.3553&0.15790&-130.3560&0.3553 \\
        \cmidrule{3-10}
        & & $(0.8,0.2)$ & 96.8678&0.5241&-131.0290&0.3498&0.00144&-131.0290&0.3498 \\
        \cmidrule{3-10}
        & & $(0.5,0.5)$ & 99.5821&0.5045&-135.7740&0.3396&0.00000&-135.7740&0.3396 \\
        \cmidrule{3-10}
        & & $(0.2,0.8)$ & 108.7352&0.4917&-149.3090&0.3332&0.00000&-149.3090&0.3332 \\
        \cmidrule{3-10}
        & & $(0,1)$ & 133.3052&0.4762&-183.8160&0.3238&0.00000&-183.8160&0.3238 \\
        \cmidrule{1-10}
        \multirow{6}{*}[-1em]{3} & \multirow{6}{*}[-1em]{0.001349} & $(1,0)$ & 96.8866&0.5300&-131.0510&0.3544&0.00130&-131.0510&0.3544 \\
        \cmidrule{3-10}
        & & $(0.9,0.1)$ & 96.7519&0.5298&-131.0510&0.3545&0.00126&-131.0510&0.3545 \\
        \cmidrule{3-10}
        & & $(0.8,0.2)$ & 96.8852&0.5240&-131.0560&0.3497&0.00106&-131.0560&0.3497 \\
        \cmidrule{3-10}
        & & $(0.5,0.5)$ & 99.6087&0.5049&-135.9450&0.3399&0.00000&-135.9450&0.3399 \\
        \cmidrule{3-10}
        & & $(0.2,0.8)$ & 108.7883&0.4916&-149.3880&0.3331&0.00000&-149.3880&0.3331 \\
        \cmidrule{3-10}
        & & $(0,1)$ & 133.9731&0.4755&-184.5450&0.3229&0.00000&-184.5450&0.3229 \\
        \bottomrule
        \label{tbl:L_Data}
    \end{tabular}
\end{table*}

\section{Discussions}
\label{disu}
In this section we closely scrutinize the two numerical examples for comparison, verification, and insights.

Visual inspection and analysis of topology optimized designs (i.e., identifying and comparing their differences) is largely an untouched topic in TO research, which, in our opinion, is curious because topology is all about geometry and ``appearance'' of structures. Without such tools, the best effort is to compare those results qualitatively. To make their differences more pronounced, readers may render them into short animations, which would reveal much more than human eyes can perceive using only static images. There would be subtle material re-distributions (i.e., among results under the same tuple $(\epsilon,1-\epsilon)$), as well as thickening or thinning of certain features, which could be almost undetectable by comparing static images. The removal or addition of features is easier to spot among results. The most distinct results correspond to $\epsilon=1$ and $\epsilon=0$, which is understandable because they are the extreme cases. The results between them are sort of transitions from one bounding value to the other.

A number of trends can be observed from both the optimized designs and their corresponding numerical results:
\begin{enumerate}
\item The MCS-based $P_f$ is always smaller than the expected $P_f$, which confirms that the designs have achieved the desired reliability level.

\item Decreasing $\epsilon$, or increasing the weight on standard deviation in the robust objective, makes the MCS-based $P_f$ smaller until reaching $0$. We hypothesize that there are two classes of solutions depending on the weight: one on the constraint boundary and the other inside the feasible region of the optimization problem. In the first class, the MCS-based $P_f$ is close to the expected $P_f$: decreasing rate of the MCS-based $P_f$ is very slow for certain range of $\epsilon$ (i.e., $\epsilon=\{1,0.9,0.8\}$ in Table \ref{tbl:C_Data}, and $\epsilon=\{1,0.9\}$ in Table \ref{tbl:L_Data}). In the second class, the solutions are away from the constraint boundary but still inside the feasible region: the decreasing rate accelerates and the MCS-based $P_f$ eventually becomes 0 (i.e., $\epsilon \geq 0.5$ in Table \ref{tbl:C_Data}, and $\epsilon \geq 0.8$ in Table \ref{tbl:L_Data}).

\item From our understanding of robust optimization, increasing the weight on standard deviation in a minimization problem will decrease its value and have the opposite effect on the mean, which is actually the case as observed in Tables \ref{tbl:C_Data} and \ref{tbl:L_Data}. However, there are a couple of outliers to this trend (i.e., $(\epsilon,\beta)=(0.8,2)$ in Table \ref{tbl:C_Data}, and $\epsilon=0.9$ in Table \ref{tbl:L_Data}). An explanation can be made by investigating where the solution converges in the design space. When the solutions are away from the constraint boundary, the robust objective is more dominant in the solution: a marked increase and decrease of the mean and standard deviation, respectively, is observed. In solutions at the constraint limit, the robust objective has less in the solution.

This trend explains the decreasing tendency of the MCS-based $P_f$. Increasing the weight on standard deviation leads to decreased influence of the mean compliance, whose major proportion is contributed by the displacement of point B (Fig.~\ref{fig:C} and \ref{fig:L}) which is constrained probabilistically in (\ref{e:examples}) to be larger than the minimum allowable value $u_0$. This results in increasing the mean of point B displacement, whose changing rate is slow for solutions close to constraint boundary and much more rapid in the remaining cases (i.e., the sixth column of Tables \ref{tbl:C_Data} and \ref{tbl:L_Data}). The bigger the displacement of point B becomes, the smaller the MCS-based $P_f$ is. It is obvious that under a specific set of inputs (i.e., loading and boundary conditions, material property) mechanical capabilities of a structure (i.e., stress, strain, displacement), even if heavily optimized, are always finite. Thus, the MCS-based $P_f$ approaches 0 when the displacement of point B continues to increase.

\item The Smolyak-type sparse grid and the SRSM are both very good methods for their respective approximating targets. As shown in Tables \ref{tbl:C_Data} and \ref{tbl:L_Data}, approximately six significant digits are required to see the differences between the MCS-based and the SRSM results. The cumulative distribution functions of point B displacement from the two methods are also almost identical, so they are not included in the paper. This also makes us confident in choosing only two terms in the KL expansion$-$adding more terms would only increase computational cost without clear benefits. The same observation is applied to the Smolyak-type sparse grid and MCS-based results of the mean and standard deviation of the compliance. Because of such agreement, we decide to use only one level of the sparse grid. Multiple levels were tested in \cite{lazarov_topology_2012,zhao_robust_2015}.
\end{enumerate}

\section{Conclusions}
\label{concl}
In this paper, we have presented an efficient robust reliability-based design approach for topology optimized designs, considering random field uncertainty in material property. Our approach avoids a double-loop approach to reliability evaluation, and utilizes a number of techniques to reduce computational cost without sacrificing solution accuracy. The two case studies have demonstrated the methodology, and have shown the impact of using the robust reliability-based design approach. As shown in the case studies, the topology can converge to a (local) solution where either the reliability-based constraint or the robust objective has more influence upon the resulting topology. This difference in converged solution is determined by the relative weight of the mean versus the variance in the robust objective. This finding is significant to designers because it shows that using either exclusively a robust design approach or reliability-based design approach will not identify the optimal topology considering material property uncertainty. Only the robust reliability-based design approach is capable of identifying the optimal topology considering the influence of material property uncertainty in both the objective and the constraint. 

Explicitly considering material property uncertainty will be significant in design for additive manufacturing, where the additive process can lead to significant material property uncertainty due to temperature gradients or other sources of variation. In terms of future work, the following are recommended. Further case studies should be considered to better understand the interplay between the reliability constraint and the robust objective. Different objectives and constraints can be considered as well. Considering topology optimization as a tool for design for additive manufacturing, another significant characteristic of the material property is non-linearity. Polymeric materials used in additive manufacturing will be better modeled as hyper-elastic or visco-elastic, as opposed to linear elastic. Another need is visualization of the results. As noted in the text, it can be difficult to understand the differences in topology based upon visual inspection. A more systematic comparison of different converged results will help researchers better understand the differences in designs resulting from different problem formulations, and will lead to a better understanding of the design principles leading to robust, reliable designs.

\begin{acknowledgment}
The first author would like to thank the Vietnam Education Foundation (VEF) for their financial support through the VEF fellowship, and Prof. Krister Svanberg for his assistance on the MMA code.
\end{acknowledgment}

\bibliographystyle{asmems4}
\bibliography{refs}

\end{document}